 
\hoffset=1cm
\voffset=1cm

\magnification=1200
\hsize=11.25cm
\vsize=18cm
\parskip 0pt
\parindent=12pt
\frenchspacing

\def \inv{\mathop{\rm inv} \nolimits}
\def \maj{\mathop{\rm maj} \nolimits}
\def \stat{\mathop{\rm stat} \nolimits}
\def \des{\mathop{\rm des} \nolimits}

\def\epsilon{\varepsilon}

\def\qed{\quad\raise -2pt\hbox{\vrule\vbox to 10pt{\hrule 
width 4pt \vfill\hrule}\vrule}}

\def\halmos{\penalty 500 \hbox{\qed}\par\smallskip}

\font\eightbf=cmbx8
\font\eightrm=cmr8 \font\sixrm=cmr6 
\font\eighttt=cmtt8

\def\diag{\mathop{\rm diag}\nolimits}
\def\card{\mathop{\rm card}\nolimits}
\def\proof{{\it Proof.\enspace }}
\def\interr{\mathop{\rm ?}\nolimits}
\def\nonfleche{\mathrel{\kern1pt\not\kern-2pt
\leftarrow\kern1pt}}
\def\bfc{{\bf c}}
\long\def\proclaim #1. #2\endproclaim{\medbreak
{\bf #1.\enspace}{\sl#2}\par\medbreak}

\catcode`\@=11

\def\eightpoint{%
  \textfont0=\eightrm \scriptfont0=\sixrm 
\scriptscriptfont0=\fiverm
  \def\rm{\fam\z@\eightrm}%
   \abovedisplayskip=9pt plus 2pt minus 6pt
  \abovedisplayshortskip=0pt plus 2pt
  \belowdisplayskip=9pt plus 2pt minus 6pt
  \belowdisplayshortskip=5pt plus 2pt minus 3pt
  \smallskipamount=2pt plus 1pt minus 1pt
  \medskipamount=4pt plus 2pt minus 1pt
  \bigskipamount=9pt plus 3pt minus 3pt
  \normalbaselineskip=9pt
  \setbox\strutbox=\hbox{\vrule height7pt depth2pt width0pt}%
  \let\bigf@ntpc=\eightrm \let\smallf@ntpc=\sixrm
  \normalbaselines\rm}

\def\appeln@te{}
\def\vfootnote#1{\def\@parameter{#1}\insert
  \footins\bgroup\eightpoint
  \interlinepenalty\interfootnotelinepenalty
  \splittopskip\ht\strutbox %
  \splitmaxdepth\dp\strutbox \floatingpenalty\@MM
  \leftskip\z@skip \rightskip\z@skip
  \ifx\appeln@te\@parameter\indent \else{\noindent #1\ }\fi
  \footstrut\futurelet\next\fo@t}

\def\footnoterule{\kern-6\p@
  \hrule width 2truein \kern 5.6\p@} 

\catcode`\@=12


 \message{`lline' & `vector' macros from LaTeX}
 \catcode`@=11
\def\{{\relax\ifmmode\lbrace\else$\lbrace$\fi}
\def\}{\relax\ifmmode\rbrace\else$\rbrace$\fi}
\def\newcount{\alloc@0\count\countdef\insc@unt}
\def\newdimen{\alloc@1\dimen\dimendef\insc@unt}
\def\newwrite{\alloc@7\write\chardef\sixt@@n}

\newwrite\@unused
\newcount\@tempcnta
\newcount\@tempcntb
\newdimen\@tempdima
\newdimen\@tempdimb
\newbox\@tempboxa

\def\@spaces{\space\space\space\space}
\def\@whilenoop#1{}
\def\@whiledim#1\do #2{\ifdim #1\relax#2%
\@iwhiledim{#1\relax#2}\fi}
\def\@iwhiledim#1{\ifdim #1\let\@nextwhile=\@iwhiledim
        \else\let\@nextwhile=\@whilenoop\fi\@nextwhile{#1}}
\def\@badlinearg{\@latexerr{Bad \string\line
\space or \string\vector
   \space argument}}
\def\@latexerr#1#2{\begingroup
\edef\@tempc{#2}\expandafter\errhelp\expandafter{\@tempc}%
\def\@eha{Your command was ignored.
^^JType \space I <command> <return> \space to replace it
  with another command,^^Jor \space <return> 
\space to continue without it.}
\def\@ehb{You've lost some text. \space \@ehc}
\def\@ehc{Try typing \space <return>
  \space to proceed.^^JIf that doesn't work, 
type \space X <return> \space to
  quit.}
\def\@ehd{You're in trouble here.  \space\@ehc}

\typeout{LaTeX error. 
\space See LaTeX manual for explanation.^^J
 \space\@spaces\@spaces
\@spaces Type \space H <return> \space for
 immediate help.}\errmessage{#1}\endgroup}
\def\typeout#1{{\let\protect\string\immediate
\write\@unused{#1}}}

\font\tenln    = line10
\font\tenlnw   = linew10

\newdimen\@wholewidth
\newdimen\@halfwidth
\newdimen\unitlength 

\unitlength =1pt


\def\thinlines{\let\@linefnt\tenln \let\@circlefnt\tencirc
  \@wholewidth\fontdimen8\tenln \@halfwidth .5\@wholewidth}
\def\thicklines{\let\@linefnt\tenlnw \let\@circlefnt\tencircw
  \@wholewidth\fontdimen8\tenlnw \@halfwidth .5\@wholewidth}

\def\linethickness#1{\@wholewidth #1\relax 
\@halfwidth .5\@wholewidth}

\newif\if@negarg

\def\lline(#1,#2)#3{\@xarg #1\relax \@yarg #2\relax
\@linelen=#3\unitlength
\ifnum\@xarg =0 \@vline
  \else \ifnum\@yarg =0 \@hline \else \@sline\fi
\fi}

\def\@sline{\ifnum\@xarg< 0 \@negargtrue 
\@xarg -\@xarg \@yyarg -\@yarg
  \else \@negargfalse \@yyarg \@yarg \fi
\ifnum \@yyarg >0 \@tempcnta\@yyarg \else 
\@tempcnta -\@yyarg \fi
\ifnum\@tempcnta>6 \@badlinearg\@tempcnta0 \fi
\setbox\@linechar\hbox{\@linefnt
\@getlinechar(\@xarg,\@yyarg)}%
\ifnum \@yarg >0 \let\@upordown\raise \@clnht\z@
   \else\let\@upordown\lower \@clnht \ht\@linechar\fi
\@clnwd=\wd\@linechar
\if@negarg \hskip -\wd\@linechar 
\def\@tempa{\hskip -2\wd\@linechar}\else
     \let\@tempa\relax \fi
\@whiledim \@clnwd <\@linelen \do
  {\@upordown\@clnht\copy\@linechar
   \@tempa
   \advance\@clnht \ht\@linechar
   \advance\@clnwd \wd\@linechar}%
\advance\@clnht -\ht\@linechar
\advance\@clnwd -\wd\@linechar
\@tempdima\@linelen\advance\@tempdima -\@clnwd
\@tempdimb\@tempdima\advance\@tempdimb -\wd\@linechar
\if@negarg \hskip -\@tempdimb \else \hskip \@tempdimb \fi
\multiply\@tempdima \@m
\@tempcnta \@tempdima \@tempdima \wd\@linechar 
\divide\@tempcnta \@tempdima
\@tempdima \ht\@linechar \multiply\@tempdima \@tempcnta
\divide\@tempdima \@m
\advance\@clnht \@tempdima
\ifdim \@linelen <\wd\@linechar
   \hskip \wd\@linechar
  \else\@upordown\@clnht\copy\@linechar\fi}

\def\@hline{\ifnum \@xarg <0 \hskip -\@linelen \fi
\vrule height \@halfwidth depth \@halfwidth width \@linelen
\ifnum \@xarg <0 \hskip -\@linelen \fi}

\def\@getlinechar(#1,#2){\@tempcnta#1\relax\multiply
\@tempcnta 8
\advance\@tempcnta -9 \ifnum #2>0 
\advance\@tempcnta #2\relax\else
\advance\@tempcnta -#2\relax\advance\@tempcnta 64 \fi
\char\@tempcnta}

\def\vector(#1,#2)#3{\@xarg #1\relax \@yarg #2\relax
\@linelen=#3\unitlength
\ifnum\@xarg =0 \@vvector
  \else \ifnum\@yarg =0 \@hvector \else \@svector\fi
\fi}

\def\@hvector{\@hline\hbox to 0pt{\@linefnt
\ifnum \@xarg <0 \@getlarrow(1,0)\hss\else
    \hss\@getrarrow(1,0)\fi}}

\def\@vvector{\ifnum \@yarg <0 \@downvector \else 
\@upvector \fi}

\def\@svector{\@sline
\@tempcnta\@yarg \ifnum\@tempcnta <0 \@tempcnta
=-\@tempcnta\fi
\ifnum\@tempcnta <5
  \hskip -\wd\@linechar
  \@upordown\@clnht \hbox{\@linefnt  \if@negarg
  \@getlarrow(\@xarg,\@yyarg) \else 
\@getrarrow(\@xarg,\@yyarg) \fi}%
\else\@badlinearg\fi}

\def\@getlarrow(#1,#2){\ifnum #2 =\z@ \@tempcnta='33\else
\@tempcnta=#1\relax\multiply\@tempcnta \sixt@@n 
\advance\@tempcnta
-9 \@tempcntb=#2\relax\multiply\@tempcntb \tw@
\ifnum \@tempcntb >0 \advance\@tempcnta \@tempcntb\relax
\else\advance\@tempcnta -\@tempcntb\advance\@tempcnta 64
\fi\fi\char\@tempcnta}

\def\@getrarrow(#1,#2){\@tempcntb=#2\relax
\ifnum\@tempcntb < 0 \@tempcntb=-\@tempcntb\relax\fi
\ifcase \@tempcntb\relax \@tempcnta='55 \or
\ifnum #1<3 \@tempcnta=#1\relax\multiply\@tempcnta
24 \advance\@tempcnta -6 \else \ifnum #1=3 \@tempcnta=49
\else\@tempcnta=58 \fi\fi\or
\ifnum #1<3 \@tempcnta=#1\relax\multiply\@tempcnta
24 \advance\@tempcnta -3 \else \@tempcnta=51\fi\or
\@tempcnta=#1\relax\multiply\@tempcnta
\sixt@@n \advance\@tempcnta -\tw@ \else
\@tempcnta=#1\relax\multiply\@tempcnta
\sixt@@n \advance\@tempcnta 7 \fi\ifnum #2<0 
\advance\@tempcnta 64 \fi
\char\@tempcnta}

\def\@vline{\ifnum \@yarg <0 \@downline \else \@upline\fi}

\def\@upline{\hbox to \z@{\hskip -\@halfwidth \vrule
  width \@wholewidth height \@linelen depth \z@\hss}}

\def\@downline{\hbox to \z@{\hskip -\@halfwidth \vrule
  width \@wholewidth height \z@ depth \@linelen \hss}}

\def\@upvector{\@upline\setbox\@tempboxa
\hbox{\@linefnt\char'66}\raise
     \@linelen \hbox to\z@{\lower 
\ht\@tempboxa\box\@tempboxa\hss}}

\def\@downvector{\@downline\lower \@linelen
      \hbox to \z@{\@linefnt\char'77\hss}}

\thinlines

\newcount\@xarg
\newcount\@yarg
\newcount\@yyarg
\newcount\@multicnt
\newdimen\@xdim
\newdimen\@ydim
\newbox\@linechar
\newdimen\@linelen
\newdimen\@clnwd
\newdimen\@clnht
\newdimen\@dashdim
\newbox\@dashbox
\newcount\@dashcnt
 \catcode`@=12

\newbox\tbox
\newbox\tboxa

\def\leftzer#1{\setbox\tbox=\hbox to 0pt{#1\hss}%
     \ht\tbox=0pt \dp\tbox=0pt \box\tbox}

\def\rightzer#1{\setbox\tbox=\hbox to 0pt{\hss #1}%
     \ht\tbox=0pt \dp\tbox=0pt \box\tbox}

\def\centerzer#1{\setbox\tbox=\hbox to 0pt{\hss #1\hss}%
     \ht\tbox=0pt \dp\tbox=0pt \box\tbox}

\def\leftput(#1,#2)#3{\setbox\tboxa=\hbox{%
    \kern #1\unitlength
    \raise #2\unitlength\hbox{\leftzer{#3}}}%
    \ht\tboxa=0pt \wd\tboxa=0pt \dp\tboxa=0pt\box\tboxa}

\def\rightput(#1,#2)#3{\setbox\tboxa=\hbox{%
    \kern #1\unitlength
    \raise #2\unitlength\hbox{\rightzer{#3}}}%
    \ht\tboxa=0pt \wd\tboxa=0pt \dp\tboxa=0pt\box\tboxa}

\def\centerput(#1,#2)#3{\setbox\tboxa=\hbox{%
    \kern #1\unitlength
    \raise #2\unitlength\hbox{\centerzer{#3}}}%
    \ht\tboxa=0pt \wd\tboxa=0pt \dp\tboxa=0pt\box\tboxa}

\unitlength=1mm

\def\cput(#1,#2)#3{\noalign{\nointerlineskip
\centerput(#1,#2){#3}
\nointerlineskip}}

\def\segment(#1,#2)\dir(#3,#4)\long#5{%
\leftput(#1,#2){\lline(#3,#4){#5}}}


\overfullrule=0pt


\vglue 1.5cm
\centerline{\bf GRAPHICAL MAJOR INDICES}
\bigskip
\centerline{\it
Dominique FOATA\/\footnote{$^1$}{D\'epartement de
math\'ematique, Universit\'e Louis Pasteur, 7, rue
Ren\'e-Descartes,\hfil\break F-67084 Strasbourg, France.
{(\eighttt foata@math.u-strasbg.fr)}. Supported in part by
the E.E.C. program on  Algebraic Combinatorics.}
and  Doron ZEILBERGER\/\footnote{$^2$}{\eightrm  \raggedright
\baselineskip=9pt Department of Mathematics, Temple University,
Philadelphia, PA 19122, USA.\hfil\break
{(\eighttt zeilberg@math.temple.edu)}. 
Supported in part by the NSF.}}

\bigskip\bigskip

{\narrower\narrower\eightpoint

\noindent
{\eightbf Abstract:} 
A generalization of the classical statistics ``maj'' and ``inv''
(the major index and number of inversions) on words is introduced,
parameterized by arbitrary graphs on the underlying alphabet.
The question of characterizing those graphs that lead to equi-distributed
``inv" and ``maj" is posed and answered.

\medskip\noindent
{\eightbf R\'esum\'e:}
On introduit une g\'en\'eralisation des statistiques classiques
que sont ``maj'' et ``inv'' (l'indice majeur et le nombre
d'inversions) sur les mots, qui est param\'etris\'ee par
des graphes arbitraires sur l'alphabet sous-jacent. La question
de caract\'eriser ces graphes conduisant \`a des statistiques
``inv" et ``maj" qui soient \'equidistribu\'ees est pos\'ee et
r\'esolue.

}

\bigskip\medskip
\centerline{\bf 0. Introduction}

\medskip
Every mathematician knows what the {\it the  number  of 
inversions} of a  permutation  is,  as it features   in  the 
definition of the determinant. The number of inversions 
of a permutation of length $n$,
$$
\inv \pi = \sum_{1 \leq i \leq j \leq n} 
\chi (\pi(i) > \pi(j)),
$$
(using the classical notation $\chi(A)=1$ or~0, depending on
whether the statement~$A$ is true of false)
is a measure of how `scrambled' it is compared to the
identity permutation $[1,2,  \dots , n]$. Netto proved (and it 
is nowadays easy to see, e.g., [Kn73, p. 15]) that the generating 
function for ``the number of inversions''
$$
\sum_{\pi \in S_n} q^{\inv (\pi )} ,
$$
equals the $q$-analog of $n!$, i.e.,
$[n]!:=1(1+q)(1+q+q^2) \dots (1+q+ \dots + q^{n-1})$,
that, can be also written as  $(q)_n /(1-q)^n$, where, as
usual in $q-$theory, $(q)_n=(1-q)(1-q^2) \dots (1-q^n)$.

The  number of inversions is an example of a {\it permutation
statistic}, by which is meant a numerical attribute that
permutations possess (just like height, weight, or number
of children for humans). The utility of the generating function
according to a given statistic ``stat,"
$$
F_{\stat}(q):= \sum_{\pi \in S_n}  q^{\stat( \pi)} ,
$$
is that it contains in it all the `statistical' information
regarding ``stat." Also its derivatives evaluated at $q=1$
enable us to, successively, find the average, standard 
deviation, and higher moments of its distribution. 
Furthermore, when the generating function is `nice' it hints at
(combinatorial, algebraic and sometimes analytic) structures.

MacMahon [Mac15, p.~135] was the first to introduce another such
statistic, that he called `the greater index', but that is
nowadays called the `major index' and denoted by ``maj." In fact,
he defined that statistic not only for permutations but for
arbitrary {\it words} with possible repetitions of letters. He did
also the same for ``inv."  If $X$ is a totally ordered alphabet,
and if $w=x_1x_2\ldots x_m$ is a word with letters in~$X$,
those two statistics are defined by
$$\eqalign{ 
\maj w&=\sum_{i=1}^{m-1} i\,\chi(x_i>x_{i+1}),\cr  
\inv w&=\sum_{1\le i<j\le m} \chi(x_i>x_j).\cr} 
$$

To restate MacMahon's result we will take the alphabet
$X$ as the linear set $[\,r\,]=\{1,2,\ldots, r\}$ $(r\ge 1)$. 
Let $\bfc=(c(1),c(2),\ldots,c(r))$ be a sequence of $r$
non-negative integers and let $v$ be the (non-decreasing) word
$v=1^{c(1)}2^{c(2)}\ldots r^{c(r)}$. We will denote by $R(v)$ 
(or by $R(\bfc)$ if there is no ambiguity) the class of all
rearrangements of the word~$v$, i.e., the class of all words
containing exactly $c(i)$ occurrences of the letter~$i$ for all
$i=1,\ldots,r$. Then MacMahon [Mac13] (see also [Mac78]) proved
that for each integer~$k$ there are as many words $w\in
R(\bfc)$ such that $\maj w=k$, as there are words $w'\in R(\bfc)$
such that $\inv w'=k$. In other words, the statistics ``maj" and
``inv" are equidistributed on each rearrangement class. 

It is well known, and easy to see, that the number of words in
$R(\bfc)$ is the multinomial coefficient:
$${c(1)+c(2)+\cdots+c(r)\choose c(1),c(2),\ldots\,,c(r)}
={(c(1)+c(2)+\cdots+c(r))!\over c(1)!\,c(2)!\, \ldots\,c(r)!}.
$$  
MacMahon's proof [Mac13, Mac78] (see also [Kn73, p. 17], [An76,
chap.~3]) of the forementioned result was to show that
the generating functions for ``inv"  and ``maj", over the class
$R(v)$, i.e., ${\sum\limits_w} q^{\inv w}$ and
${\sum\limits_w} q^{\maj w}$ (with $w$ runnning over
the class $R(\bfc)$), were {\it both} given by the
$q$-analog of the multinomial coefficient: 
$${c(1)+c(2)+\cdots+c(r)\brack c(1),c(2),\ldots\,,c(r)}
={(q)_{c(1)+c(2)+\cdots+c(r)}\over 
(q)_{c(1}(q)(_{c(2)}\ldots (q)_({c(r)}}.
$$    
 
The natural question of finding a bijection that
sends each permutation to another one in such a way that the
major index of the image equals the number of inversions of
the original, has been answered by the first author [Fo68], 
and since `canonized ' in {\it the book} ([Kn73], ex. 5.1.1.19).
 
In this paper we introduce a natural generalization of both 
``inv" and ``maj," parameterized by a general directed graph. 
A {\it directed graph} on $X$ is any subset $U$ of the
Cartesian product $X\times X=
\{(x,y) \,\vert\, 1 \leq x \leq r, 1 \leq y \leq r \}$.
Of course there are altogether $2^{r^2}$ directed graphs.

For each such directed graph $U$ let's associate the following
statistics defined on each word
$w=x_1x_2\ldots x_m$ by
$$\eqalign{ 
\maj'_U w&=\sum_{i=1}^{m-1} 
i\,\chi\bigl((x_i,x_{i+1})\in U\bigr),\cr   
\inv'_U w&=\sum_{1\le i<j\le m} 
\chi\bigl((x_i,x_{i+1})\in U\bigr).\cr}\leqno(0.1) 
$$
Further in the paper other statistics ``$\maj_U$" and ``$\inv_U$"
(without any primes) will be introduced.
 
The purpose of this paper is to characterize the
directed graphs $U$ that posses the `Mahonian property' of ``inv"
and ``maj" having the same generating function. We first need the
following definition. 

\medskip
{\bf Definition}.
An {\it ordered bipartition} of $X$ is a sequence
$(B_1,B_2,\ldots, B_k)$ of non-empty disjoint subsets of~$X$, of
union~$X$, together with a sequence 
$(\beta_1,\beta_2,\ldots,\beta_k)$ of elements equal to 0 or~1. 
If $\beta_l=1$ (resp.~0), we say that the subset $B_l$ is 
{\it underlined} (resp. {\it non-underlined}). For the sake of
convenience, we also say that the {\it subscript}~$l$ or 
{\it each element} of~$B_l$ is {\it underlined} (resp. 
{\it non-underlined}).
 
A relation $U$ on $X\times X$ is said to be {\it bipartitional},
if there exists an ordered bipartition 
$
((B_1,B_2,\ldots, B_k)
,\ (\beta_1,\beta_2,\ldots,\beta_k))
$
such that
$(x,y)\in U$ iff either $x\in B_l$, $y\in B_{l'}$ and $l<l'$, i.e.,
if the block containing~$x$ is to the left of the block
containing~$y$, or $x$ and $y$ belong to the same block~$B_l$
and $B_l$ is underlined.

\medskip
As proved by Han [Han95], a bipartitional relation~$U$ can also
be characterized by the following two relations
$$
\eqalign{(x,y)\in U,\ (y,z)\in U&\Rightarrow (x,z)\in U\,;\cr
(x,y)\not\in U,\ (z,y)\in U&\Rightarrow (z,x)\in U.\cr} 
$$
Some particular bipartitional relations are worth being noticed.

1) $U=\{(x,y)\,|\,x>y\}$ that corresponds to the ordered
bipartition $(\{r\},\ldots,\{2\},\{1\})$ ; in this case
$\inv_U=\inv$ and $\maj_U=\maj$;

2) $U=\{(x,y)\,|\,x\ge y\}$ that is associated with the ordered
bipartition
$(\underline{\{ r\}},\ldots,\underline{\{ 2\}},
\underline{\{ 1\}})$,
where all the blocks are underlined; the inversions
and descents involved in the statistics ``$\inv_U$" and
``$\maj_U$'' also include all the pairs $(x,x)$;

3) $U=\emptyset$ which is associated with the one-block ordered
bipartition $(\{1,2,\ldots,r\})$; the statistics ``$\inv_U$" and
``$\maj_U$" are identically zero;

4) $U=X\times X$ which is associated with the 
one-underlined-block ordered bipartition 
$(\underline{\{1, 2,\ldots, r\}})$; in this case $\inv_U w=\maj_U
w=m(m-1)/2$ for each word $w$ of length~$m$;

5) $U$ which is associated with an ordered bipartition all the
blocks of which are singletons; such relations have been
considered by Clarke and Foata [ClFo94, ClFo95a, ClFo95b] who
also introduced the statistic ``$\maj_k$" which is immediately
related with the statistic ``$\maj_U$" further defined.

\medskip
A bipartitional relation 
$U=((B_1,B_2,\ldots, B_k),\ (\beta_1,\beta_2,\ldots,\beta_k))$
can also be visualized as follows:
rearrange the elements of~$X$ in a row in such a way that the
elements of $B_1$ come first, in any order, then the elements
of~$B_2$, etc. Then $U$ will consist of all the block products
$B_l\times B_{l'}$ with $l<l'$, as well as the block
product $B_l\times B_l$ whenever $B_l$ is underlined.

\medskip
In Figure 1, for instance, the underlying ordered bipartition
consists of four blocks $(B_1,B_2,B_3,B_4)$ with $B_1$, $B_4$
underlined.

\midinsert
\def\segment(#1,#2)\dir(#3,#4)\long#5{%
\leftput(#1,#2){\lline(#3,#4){#5}}}

\newbox\partitionalA

\setbox\partitionalA=\vbox{\offinterlineskip 
\segment(0,0)\dir(0,-1)\long{40}
\segment(0,0)\dir(1,0)\long{40}
\segment(40,0)\dir(0,-1)\long{40}
\segment(0,-40)\dir(1,0)\long{40}
\segment(5,0)\dir(0,-1)\long{40}
\segment(0,-35)\dir(1,0)\long{40}
\segment(0,-15)\dir(1,0)\long{40}
\segment(25,0)\dir(0,-1)\long{40}
\segment(0,-20)\dir(1,0)\long{40}
\segment(20,0)\dir(0,-1)\long{40}
\centerput(3,-39){$U$}
\centerput(3,-29){$U$}
\centerput(3,-19){$U$}
\centerput(3,-10){$U$}
\centerput(12,-10){$U$}
\centerput(22,-10){$U$}
\centerput(32,-10){$U$}
\centerput(12,-19){$U$}
\centerput(-5,-10){$B_4$}
\centerput(-5,-19){$B_3$}
\centerput(-5,-29){$B_2$}
\centerput(-5,-39){$B_1$}
\centerput(3,-45){$B_1$}
\centerput(12,-45){$B_2$}
\centerput(22,-45){$B_3$}
\centerput(32,-45){$B_4$}
}

\centerline{\box\partitionalA\hskip
37mm} 
\kern 135pt
\centerline{Fig. 1}

\endinsert

Our first result is the following.
 
\proclaim Theorem 1. The statistics ``$\inv'_U$" and ``$\maj'_U$"
are equidistributed on each rearrangement class, if and only
if the relation~$U$ is bipartitional.
\endproclaim

We first prove the `easy' part, which as usual is the `if'
part. Three proofs will be given. The first manipulative, the
second combinatorial \`a la MacMahon, the third bijective, as
people say to-day. All this is derived in sections 3, 4 and~5,
respectively. Section~6 contains the proof of the `only if' part.

Now if $U$ is a bipartitional relation on~$X$, two other
statistics ``$\inv_U$" and ``$\maj_U$" may be defined, that
also reduce to ``inv" and ``maj" when
$U=(\{r\},\ldots,\{2\},\{1\})$. 
Let $|w|_{-}$ denote the number of underlined letters
in the word $w=x_1x_2\ldots x_m$. Then define 
$$\eqalign{\
\maj_U w&=\sum_{i=1}^{m-1} 
i\,\chi\bigl((x_i,x_{i+1})\in U)\bigr)
+m\,\chi(x_m\ {\rm is\ underlined}),\cr   
\inv_U w&=\sum_{1\le i<j\le m} 
\chi\bigl((x_i,x_j)\in U)\bigr)+| w|_{-}\,.
\cr} \leqno(0.2)
$$
In other words, $\maj_U$ is equal to $\maj'_U w$ plus the {\it
length} of the word whenever the last letter is underlined, while
$\inv_U w$ is equal to $\inv'_U w$ plus the {\it number of
underlined letters} in~$w$.

We say that a bipartitional relation~$U$ is {\it
compatible}, if all its underlined blocks are {\it on the left} of
its non-underlined ones, or, with the above notations, if the
sequence $(\beta_1,\beta_2,\ldots,\beta_k)$ is of the form
$(1,1,\ldots,1,0,0,\ldots,0)$.

We next prove the theorem.

\proclaim Theorem 2. Let $U$ be a bipartitional relation on
$X$. Then ``$\maj_U$" and ``\/$\inv_U$" are equidistributed on
each rearrangement class, if and only if $U$ is compatible.
\endproclaim

As we shall see, the notion of compatibility is crucial. It 
relates with an analogous notion introduced in Clarke and Foata
({\it op. cit.}) for dealing with the number of {\it excedances}
and the {\it Denert statistic}.
If $U$ is non-compatible, ``$\maj_U$" and ``$\inv_U$" are not
even equidistributed on a class of two elements. For example,
let $X=\{1,2\}$ and let $U$ be the
(non-compatible) bipartitional relation associated with the
ordered bipartition $(\{1\},\underline{\{2\}})$. Then
$\inv_U 1,2=2$, $\inv_U 2,1=1$, while $\maj_U 1,2=3$,
$\maj_U 2,1=0$. Actually, that simple example is the core of the
proof of the `only if' part of Theorem~2 (see section~7).

Let $U$ be an ordered bipartition.
Parallel to the definition of ``$\maj'_U$" and ``$\maj_U$" we can
also introduce two kinds of {\it $U$-descents}. Let
$w=x_1x_2\ldots x_m$ be a word; we say that there is a
$U$-{\it descent of the first kind} at~$i$ in~$w$, if $1\le
i\le m-1$ and $(x,x_{i+1})\in U$, and a $U$-{\it descent of the
second kind} at~$i$ in~$w$, if $1\le i\le m-1$ and
$(x,x_{i+1})\in U$ or $i=m$ and $x_m$ is underlined. Denote by
$\des'_U w$ (resp. $\des_U w$) the number of those
$U$-descents of the first kind (resp. of the second kind). In
section~4 we derive an expression
for the generating function for each rearrangement class
$R(\bfc)$ by the pair of statistics $(\des'_U,\maj'_U)$. 

Section 7 contains the calculation of the generating function
of $R(\bfc)$ by the pair $(\des_U,\maj_U)$ and also the proof of
Theorem~2. A bijective proof of the latter Theorem appears in
section~8.

\bigskip
\centerline{\bf 1. Enumerating bipartitional relations}

\medskip
For each $r\ge 1$ let $b_r'$ (resp. $b_r$) be the number of
bipartitional relations (resp. compatible bipartitional
relations) on a set of cardinality~$r$. Also let
$b'_0=b_0=1$. The exponential generating functions for both
sequences $(b'_r)$ and $(b_r)$ are easily derived and, using
{\eightrm MAPLE}, their first values calculated. Denote by $S(r,k)$
$(1\le r\le k)$ the sequence of the Stirling numbers of the
second kind (see, e.g., [Co70, vol.~2, p.~40]).

\proclaim Proposition {1.1}. We have the formulas
$$\leqalignno{
b'_r&=\sum_{k=1}^r S(r,k)k!\,2^k\,;&(1.1)\cr
b_r&=\sum_{k=1}^r S(r,k)\,(k+1)! \,;\kern4.5cm&(1.2)\cr
\sum_{r\ge 0} b'_r {u^r\over r!}
&={1\over 3-2e^u}&(1.3)\cr
\noalign{\hbox{\qquad\qquad${} =
\displaystyle 1+2u+10{u^2\over 2!}+74{u^3\over
3!} +730{u^4\over 4!}+9002{u^5\over 5!}
+133210{u^6\over 6!}+\cdots$}}
\sum_{r\ge 0} b_r{u_r\over r!} 
&={1\over (2-e^u)^2}&(1.4)\cr
\noalign{\hbox{\qquad\qquad${} =
\displaystyle 1+2u+8{u^2\over 2!}+66{u^3\over 3!}
+308{u^4\over 4!}+2612{u^5\over 5!}
+25988{u^6\over 6!}+\cdots$}}}
$$ 
\endproclaim

\proof
Formulas (1.1) and (1.2) follow immediately from the
combinatorial definition of the Stirling numbers. Accordingly,
we can easily derive (1.3) and (1.4) from the ``vertical"
exponential generating function for the Stirling numbers.
A more direct and conceptual proof consists of making use of
the partitional complex approach [Fo74] (or invoking the theory
of species dear to our qu\'eb\'ecois friends [Be94]). This goes
as follows.

Suppose that for each $r\ge 1$ there are two blocks
of size~$r$, say, the underlined $\underline{[\, r\,]}$ and the
non-underlined block $[\,r\,]$. The exponential generating
function for those two kinds of blocks is
$$
G=2{u^1\over 1!}+2{u^2\over 2!}
+\cdots+2{u^n\over n!}+\cdots=2e^u-2.
$$
Hence the expansion of $(1-G)^{-1}$ will be the generating
function for the {\it ordered} sequences of blocks, some of
them being underlined and the others being non-underlined, i.e.,
for the {\it ordered bipartitions}. Furthermore,
$(1-G)^{-1}=1/(3-2e^u)$.

\goodbreak
For the compatible bipartitional partitions there are again two
kinds of blocks, but this time the underlined blocks must lie to
the left of the non-underlined ones. The exponential
generating functions for the underlined blocks and for the
non-underlined blocks are the same:
$$
H={u^1\over 1!}+{u^2\over 2!}+\cdots+{u^n\over
n!}+\cdots=(e^u-1).
$$
Hence the expansion of $(1-H)^{-1}(1-H)^{-1}$ will be the
generating function for the ordered sequences of blocks, the
leftmost ones being underlined, the rightmost ones being
non-underlined, so that
$$
\sum_{r\ge 0} b_r {u^r\over r!}
=\Bigl({1\over (1-(e^u-1))} \Bigr)^2={1\over (2-e^u)^2}.\qed
$$

The sequences $(b'_r)$ and $(b_r)$ do not appear (yet?) in
the Sloane integral sequence basis [Sl94]. However our young
colleague Jiang Zeng drew our attention to the paper by Knuth
[Kn92] who himself pointed out that the generating function
$(2-e^u)^{-1}$ already appeared in Cayley ({\sl Collected Math.
Papers}, vol.~4, p.~112-115) for enumerating a special class of
trees. According to Knuth the coefficients of the Taylor
expansion of $(2-e^u)^{-1}$ count the {\it preferential
arrangements of $n$ objects}.

\bigskip
\centerline{\bf 2. Notations and first analytic results}

\medskip
We make use of the usual notations:
$(a;q)_n$ and $(a;q)_{\infty}$ for the $q$-ascending factorials:
$$
\eqalign{
(a;q)_n&=\cases{1,&if $n=0$;\cr
(1-a)(1-aq)\ldots (1-aq^{n-1}),&if $n\ge 1$;\cr}\cr
(a;q)_{\infty}&=\lim\nolimits_n(a;q)_n
=\prod\limits_{n\ge 0} (1-aq^n).\cr}
$$
In particular, $(q)_n=(q;q)_n$ and $(q)_\infty=(q;q)_\infty$.

Recall the $q$-binomial theorem (see [An76
p.~15] or [GaRa90, \S\kern2pt 1.3]) that states
$$\leqalignno{
\sum_{n\ge0}{(a;q)_n\over (q;q)_n}u^n
&={(au;q)_\infty\over (u;q)_\infty},&(2.1)\cr
\noalign{\hbox{together with the two $q$-exponential
identities}} 
e_q(u)&=\sum_{n\ge 0} {u^n\over
(q;q)_n}={1\over (u;q)_\infty};&(2.2)\cr
E_q(u)&=\sum_{n\ge 0} {q^{n\choose 2} u^n\over
(q;q)_n}=(-u;q)_\infty.&(2.3)\cr
}
$$
The $q$-binomial theorem provides the five
expansions (see [An76, p.~15])
$$\leqalignno{
\sum_{n\ge 0} {s+n\brack n} u^n 
&={1\over (u;q)_{s+1}};&(2.4)\cr
\sum_{n\ge 0} {s\brack n} q^{n+1\choose 2}u^n
&=(-qu;q)_s.&(2.5)\cr
}
$$
$$\leqalignno{\noalign{\vskip-12pt}
{s+n\brack n}&=\sum_{s\ge a_1\ge\cdots \ge a_n\ge 0}
q^{a_1+\cdots+a_n}
=\sum_{n\ge a_1\ge\cdots \ge a_s\ge 0}
q^{a_1+\cdots+a_s}\,;&(2.6)\cr
\qquad q^{n\choose 2}{s+1\brack n}
&=\sum_{s\ge a_1>\cdots>a_n\ge 0}
q^{a_1+\cdots +a_n}.&(2.7)\cr
{1\over (u;q)_{s+1}}
&=\sum_{n\ge 0} u^n \sum_{n\ge a_1\ge\cdots \ge a_s\ge 0}
q^{a_1+\cdots+a_s}\,;&(2.8)\cr
}
$$
where the $a_i$'s are non-negative integers.

The ordered bipartition $((B_1,B_2,\ldots,B_k),(\beta_1,\beta_2,
\ldots, \beta_k))$ will be kept fixed throughout this section.
Let $U$ be the bipartitional relation on
$X\times X$ associated with it. Next consider a sequence
${\bf c}=(c(1),c(2),\ldots,c(r))$ of non-negative integers; as
before, let $v=1^{c(1)}2^{c(2)}\ldots r^{c(r)}$ and denote by
$R(v)$ (or by $R({\bf c})$) the class of all rearrangements of
the word~$v$. If the block $B_l$ consists of the integers
$i_1,i_2,\ldots,i_h$ written in increasing order (with respect to
the usual linear order of $X=[\,r\,]$) and if $u_1,u_2,\ldots,u_r$
are $r$ commuting variables, it will be convenient to write 
$$
\leqalignno{c(B_l)&\hbox{ for the {\it sequence}\quad}
c(i_1),c(i_2),\ldots, c(i_h)\,;&(2.9)\cr
\textstyle m_l=\sum c(B_l)&\hbox{ for the {\it sum}\quad}
c(i_1)+c(i_2)+\cdots+c(i_h)\,;\cr
\left|\bfc\right|&\hbox{ for the {\it sum} }
c(1)+\cdots+c(r)\hbox{ also equal to }m_1+\cdots
+m_k\,;\cr 
{\bf c}\pm 1_i&\hbox{ for }(c(1),\ldots,c(i-1),c(i)\pm
1,c(i+1),\ldots,c(r))\,;\cr 
u(B_l)^{c(B_l)}&\hbox{ for the {\it
monomial}\quad} u_{i_1}^{c(i_1)} u_{i_2}^{c(i_2)}\ldots
u_{i_h}^{c(i_h)}\,;\cr 
\textstyle \sum u(B_l)&\hbox{ for the
{\it sum}\quad}  u_{i_1}+u_{i_2}+\cdots +u_{i_h}.\kern 2cm\cr
{\bf u}^{\bf c}&\hbox{ for the {\it monomial} }
u_1^{c_1}u_2^{c_2}\ldots u_r^{c_r}.\cr
}
$$
In particular $\displaystyle \smash{m_l\choose c(B_l)}$ will
denote the multinomial coefficient
$$\displaystyle {c(i_1)+c(i_2)+\cdots+c(i_h)\choose
c(i_1),c(i_2),\ldots\,,c(i_h)}.
$$
Let $A'_U(q;{\bf c}):=\sum\limits_w q^{\inv'_U w}\ 
(w\in R({\bf c}))$
denote the generating function for the class $R({\bf c})$ by the
statistic ``$\inv_U'$."

\proclaim Proposition {2.1}. With the above notations $(2.9)$ the
following formulas hold 
$$
\displaylines{\noalign{\vskip-10pt}
\rlap{\rm (2.10)}\hfill
 A'_U(q;{\bf c})={|\bfc| \brack m_1,\ldots, m_k}
\prod_{l=1}^k {m_l \choose c(B_l)} 
q^{\beta_l{m_l\choose 2}};\hfill\cr
\rlap{\rm (2.11)}\hfill
\eqalign{
\sum_{\bf c} {A'_U(q;{\bf c})\over
(q)_{\left |{\bf c}\right|}} {\bf u}^{\bf c}
&=\prod_{l;\beta_l=0} e_q(\scriptstyle\sum \textstyle u(B_l))
\times \displaystyle \prod_{l;\beta_l=1}
E_q(\scriptstyle\sum \textstyle u(B_l))\hfill\cr 
&={\prod\limits_{l;\beta_l=1} (-\sum u(B_l);q)_\infty
\over \prod\limits_{l;\beta_l=0} (\sum u(B_l);q)_\infty}.\cr
}\hfill\cr
} 
$$
In the second formula {\bf c} runs over all sequences
$(c(1),\ldots,c(r))$ with $c(1)\ge 0$, \dots~, $c(r)\ge 0$.
\endproclaim

\proof
Formula (2.10) follows from the well-known generating function
in the ordinary ``inv" case. The $q$-multinomial coefficient  is
the generating function for the class of words having exactly
$m_1$ letters equal to~1, \dots~, $m_k$
letters equal to~$k$ by ``inv." Such a word gives rise to exactly
$\prod\limits_l {m_l\choose c(B_l)}$ words
in  $R({\bf c})$. Finally, the letters belonging to each
non-underlined block provide no further
$U$-inversions, while the letters in an underlined block $B_l$ 
$(\beta_l=1)$ bring 
$\smash{m_l\choose 2}$ 
extra $U$-inversions when they are compared between
themselves.

To derive (2.11) we have to make use of the traditional
$q$-calculus. First rewrite (2.10) as
$$\displaylines{\noalign{\vskip -12pt}
\rlap{(2.12)}\hfill
{A'_U(q;{\bf c})\over (q)_{\left |{\bf c}\right|}} 
=\prod_{l=1}^k 
{q^{\beta_l{m_l\choose 2}}
\over (q)_{m_l}}
{m_l \choose c(B_l)}.\hfill\cr
\noalign{\hbox{Then the left-hand side of (2.10) is equal to}}
\kern.5cm
\sum_{c(B_1)}\!\!\!\cdots\!\!\!\sum_{c(B_k)}
{A'_U(q;{\bf c})\over (q)_{\left |{\bf c}\right|}} 
u(B_1)^{c(B_1)}\ldots u(B_k)^{c(B_k)}\hfill\cr
\noalign{\vskip-10pt }
\kern3cm{}=\prod_{l=1}^k \sum_{c(B_l)}
{q^{\beta_l{m_l\choose 2}}
\over (q)_{m_l}}
{m_l \choose c(B_l)}u(B_l)^{c(B_l)}\hfill\cr
\kern3cm{}=\prod_{l=1}^k \sum_{d(l)\ge 0}
{q^{\beta_l{d(l)\choose 2}} \over (q)_{d(l)}}
\sum_{m_l=d(l)}
{m_l \choose c(B_l)}u(B_l)^{c(B_l)}\hfill\cr
\kern3cm
{}=\prod_{l=1}^k\sum_{d(l)\ge 0}
{q^{\beta_l{d(l)\choose 2}} \over (q)_{d(l)}}
\Bigl(\textstyle\sum u(B_l)\Bigr)^{d(l)},
\hfill\cr
}
$$
which is the right-hand side of (2.11) by using (2.2) and
(2.3).\halmos

\goodbreak
\bigskip
\centerline{\bf 3. An `Essentially Verification' 
Manipulative Proof}
\centerline{\bf of the `If' Part of Theorem 1}

\medskip
The generating function according to ``$\maj'_U$"  does not seem
to be directly derivable from the classical, MacMahon formula.
Later, we will show that the combinatorial proofs easily carry
over, but here we will show a manipulative proof. We will prove
the stronger result that the subsets of words with a prescribed
last letter have the Mahonian property. 

\medskip
Keep the same notations as in Proposition~2.1. In particular, let
$|{\bf c}|=c(1)+\cdots+c(r)= m_1+\cdots+m_l $ be the length of
the words in the class~$R({\bf c})$. Also define, for
each letter $i\in X=[\,r\,]$.  
$$\eqalign{
A'_U(q;{\bf c};i)&:=
\sum_w q^{\inv'_U w}\qquad
(w\in R({\bf c}),\ w\ {\rm ends\ with}\ i)\,;\cr
A''_U(q;{\bf c};i)&:=
\sum_w q^{\maj'_U w}\qquad
(w\in R({\bf c}),\ w\ {\rm ends\ with}\ i).\cr}
$$
It is easy to derive a formula for $A'_U(q;{\bf c};i)$, in terms of 
$A'_U(q;{\bf c})$. Let $i$ belong to the block $B_l$ 
$(1\le l \le k)$. Then,
$$\leqalignno{
A'_U(q;{\bf c};i)&=
A'_U(q;{\bf c}-1_i)
q^{m_1+\cdots+m_{l-1}} ,
&(3.1)\cr
\noalign{\hbox{if $B_l$ is {\it not} underlined, and}}
A'_U(q;{\bf c};i)&=
A'_U(q;{\bf c}-1_i)
q^{m_1+\cdots+m_l-1} ,
&(3.2)\cr}
$$
if $B_l$ is underlined.

Denote by $A''_U(q;{\bf c};i)$ the corresponding quantities
for ``$\maj'_U$." By considering what letter can be
second-to-last, we get the following recurrence:
$$\displaylines{(3.3)\ 
A''_U(q;{\bf c};i)=
q^{\left |{\bf c}\right|-1}\!\!\!\!\!
\sum_{j \in B_1 \cup \cdots \cup B_{l-1}}\!\!\!
A''_U(q;{\bf c}-1_i;j)
+\!\!\!\sum_{j \in B_l \cup \cdots \cup B_k}\!\!\!
A''_U(q;{\bf c}-1_i;j)\hfill\cr
\noalign{\hbox{when $B_l$ is {\it not} underlined, and the
recurrence}}
(3.4)\ A''_U(q;{\bf c};i)=
q^{\left |{\bf c}\right|-1}\!\!\!\!\!
\sum_{j \in B_1 \cup \cdots \cup B_{l}}\!\!\!
A''_U(q;{\bf c}-1_i;j)+\!\!\!\!\!\!
\sum_{j \in B_{l+1} \cup \cdots \cup B_k}\!\!\!
A''_U(q;{\bf c}-1_i;j)\hfill\cr}
$$
when $B_l$ is underlined.

It is a completely routine matter, that we leave to the readers
(or rather to their computers) to verify that the expressions
on the right sides of (3.1) and (3.2) (using (2.9)) also
satisfy the same recurrence. It follows by induction that  for
all {\bf c} and for all $1 \leq i \leq r$, we have
 $$
A'_U(q;{\bf c};i)=A''_U(q;{\bf c};i).\leqno(3.5)
$$
By summing over $i$, we get that indeed for bipartitional
graphs $U$ the statistics ``$\inv_U$" and ``$\maj_U$" are
equidistributed.\halmos
 
\bigskip
\centerline{\bf 4. The MacMahon Verfahren}

\medskip
In this section we make use again of the same
notations as in section~2. Consider the polynomial 
$$
A'_U(t,q;{\bf c})=\sum_w t^{\des'_U w} q^{\maj'_U w}\qquad
(w\in R({\bf c}))
$$
(where $\des_U' w$ is the number of $U$-descents of the
first kind in~$w$ defined in the introduction). We first derive
the formula 
$$
{A'_U(t,q;\bfc)\over (t;q)_{|\bfc|+1}}
=\prod_{l=1}^k 
{ m_l\choose c(B_l)}
\times \sum_{s\ge 0} t^s \prod_{l;\beta_l=0}
{m_l+s\brack m_l}\times
\prod_{l;\beta_l=1} q^{m_l\choose 2} {s+1\brack m_l}.
\leqno(4.1)
$$
The previous formula is the {\it finite version} of (2.11). By
multiplying (4.1) by $(1-t)$ and making $t=1$ we recover (2.11). 
>From (4.1) we also derive the {\it factorial} generating function
for the polynomials $A'_U(t,q;\bfc)$ that reads
$$
\sum_\bfc {A_U(t,q;\bfc)\over (t;q)_{|\bfc|+1}} {\bf u}^\bfc
=\sum_{s\ge 0} t^s\, {\prod\limits_{l;\beta_l=1}
(-\sum u(B_l);q)_{s+1}\over
\prod\limits_{l;\beta_l=0}
(\sum u(B_l);q)_{s+1}},\leqno(4.2)
$$
which is the {\it finite version} of (2.11).

As done in Proposition 2.1 we can obtain (4.2) from (4.1) by a
routine calculation as follows:
$$
\displaylines{\quad
\sum_\bfc {A'_U(t,q;\bfc)\over (t;q)_{|\bfc|+1}}
=\sum_{c(B_1),\ldots,c(B_k)} 
{A'_U(t,q;\bfc)\over (t;q)_{|\bfc|+1}}
\prod_{l=1}^k u(B_l)^{c(B_l)}\hfill\cr
\qquad\qquad{}=
 \sum_{c(B_1),\ldots,c(B_k)} 
\prod_{l=1}^k u(B_l)^{c(B_l)}
{ m_l\choose c(B_l)}\hfill\cr
\kern3cm{}\times \sum_{s\ge 0} t^s \prod_{l;\beta_l=0}
{m_l+s\brack m_l}\times
\prod_{l;\beta_l=1} q^{m_l\choose 2} {s+1\brack m_l}
\hfill\cr
\qquad\qquad{}=
\sum_{s\ge 0} t^s \prod_{l;\beta_l=0}\sum_{c(B_l)}
u(B_l)^{c(B_l)}{ m_l\choose c(B_l)}{m_l+s\brack m_l}\hfill\cr
\kern3cm{}\times 
\prod_{l;\beta_l=1}\sum_{c(B_l)}u(B_l)^{c(B_l)}
{ m_l\choose c(B_l)} {s+1\brack m_l}\hfill\cr
\qquad\qquad{}=
\sum_{s\ge 0} t^s \prod_{l;\beta_l=0}
\sum_{m_l}{m_l+s\brack m_l}
\sum_{\Sigma c(B_l)=m_l} 
u(B_l)^{c(B_l)}{ m_l\choose c(B_l)}\hfill\cr 
\kern3cm{}\times 
\prod_{l;\beta_l=1}\sum_{m_l} {s+1\brack m_l}
\sum_{\Sigma c(B_l)=m_l}u(B_l)^{c(B_l)}
{ m_l\choose c(B_l)}\hfill \cr
\qquad\qquad{}=
\sum_{s\ge 0} t^s \prod_{l;\beta_l=0}
\sum_{m_l}{m_l+s\brack m_l}
\Bigl({\textstyle \sum}\,u(B_l)\Bigr)^{m_l}\hfill\cr 
\kern3cm{}\times 
\prod_{l;\beta_l=1}\sum_{m_l} {s+1\brack m_l}
\Bigl({\textstyle \sum}\,u(B_l)\Bigr)^{m_l}\hfill\cr 
\qquad\qquad{}=
\sum_{s\ge 0} t^s\, {\prod\limits_{l;\beta_l=1}
(-\sum u(B_l);q)_{s+1}\over
\prod\limits_{l;\beta_l=0}
(\sum u(B_l);q)_{s+1}},\hfill\cr
}
$$
by using the identities (2.4) and (2.5).\halmos

\medskip
Now we can prove (4.1) using the so-called ``MacMahon
Verfahren." As already noted in [Fo95, ClFo95a], the method
introduced by MacMahon [Mac13] to derive the generating function
for ``maj" is to be updated to include a second statistic, but the
principle remains the same.

Let $((B_1,\ldots,B_k)(\beta_1,\ldots,\beta_k))$ be the ordered
bipartition corresponding to the bipartitional relation~$U$
and let $w=x_1x_2\ldots x_m$ be a word of the class $R(\bfc)$,
so that $m=|\bfc|$. Denote by $w(B_l)$ be the {\it subword}
of~$w$ consisting of all the letters belonging to~$B_l$
$(l=1,\ldots, k)$. Then replace each letter belonging to~$B_l$ by
$b_l=\min B_l$ (with respect to the usual order). Call 
$\overline w=\overline x_1\overline x_2\ldots \overline x_m$
the resulting word. Clear the mapping
$$
w\mapsto \bigl(\overline w, w(B_1),\ldots,w(B_l)\bigr)\leqno(4.3) 
$$
is bijective. Moreover, $\des'_U w=\des'_U\overline w$ and
$\maj'_U w=\maj'_U\overline w$. Accordingly, the polynomial
$A'_U(t,q;\bfc)$ is divisible by 
$\prod\limits_l {m_l\choose c(B_l)}$.

For each $i=1,2,\ldots,m$ let $z_i$ denote the number of
$U$-descents (of the first kind) in the right factor $\overline
x_i \overline x_{i+1}\ldots \overline x_m$ of~$\overline w$.
Clearly, $z_1=\des'_U \overline w$ and 
$z_1+\cdots +z_m=\maj'_U \overline w$.

Now let ${\bf p}=(p_1,\ldots, p_m)$ be a sequence of $m$ integers
satisfying  $s'\ge p_1\ge p_2\ge \cdots\ge p_m\ge 0$, where $s'$
is a {\it given} integer. Form the {\it non-increasing} word
$v=y_1y_2\ldots y_m$ defined by $y_i=p_i+z_i$ $(1\le i\le m)$ 
and consider the biword
$$
\Bigl(\matrix{v\cr \overline w\cr}\Bigr)=
\Bigl(\matrix{y_1y_2\ldots y_m\cr
\overline x_1\overline x_2\ldots \overline x_m\cr}\Bigr). 
$$
Next rearrange the columns of the previous matrix in such a
way that the mutual orders of the columns with the same bottom
entries are preserved and the entire bottom row is of the form
$b_1^{m_1}b_2^{m_2}\ldots b_k^{m_k}$. We obtain the matrix
$$
\Bigl(\matrix{a_{1,1}\ldots  a_{1,m_1}\cr
b_1\  \ldots\  b_1\cr}
\matrix{\ldots\cr \ldots\cr}
\matrix{a_{k,1}\ldots  a_{k,m_k}\cr
b_k\ \ldots\ b_k\cr}\Bigr).
$$
By construction each of the $k$ words
$a_{1,1}\ldots  a_{1,m_1}$, \dots~, 
$a_{k,1}\ldots  a_{k,m_k}$ is {\it non-increasing}. Furthermore,
if $\overline x_i=\overline x_{i'}$ and $\overline x_i\in B_l$
with $l$ underlined, there is necessarily a $U$-descent within
$\overline x_i \overline x_{i+1}\ldots \overline x_{i'}$. Hence
$z_i>z_{i'}$ and and $y_i>y_{i'}$. The corresponding word
$a_{l,1}\ldots a_{l,m_l}$ will then be {\it strictly decreasing}.
Also note that
$$
a_{l,i}\le y_1=p_1+z_1\le s'+\des'_U \overline w
$$
for all $l,i$. Let then $s=s'+\des'_U \overline w$. It follows
that each of the words $a_{l,1}\ldots a_{l,m_l}$ satisfies
$$
\eqalign{s&\ge a_{l,1}\ge \cdots \ge a_{l,m_l}\ge 0,
\hbox{ if $l$ is {\it not} underlined;}\cr
s&\ge a_{l,1}>\cdots >a_{l,m_l}\ge 0,
\hbox{ if $l$ is underlined.}\cr
}\leqno(4.4)
$$
The mapping $(s',{\bf p},\overline w)\mapsto (s,(a_{l,i}))$ is a
bijection satisfying
$$
\eqalign{s&=s'+\des'_U \overline w\,;\cr
\sum_{l,i}a_{l,i}=p_1+\cdots +p_m&+z_1+\cdots +z_m
=\sum_i p_i +\maj'_U \overline w.\cr}\leqno(4.5)
$$
Now rewrite (2.8) as
$$
{1\over (t;q)_{m+1}} 
=\sum_{s\ge 0} \sum_{s\ge p_1\ge\cdots\ge p_m\ge 0} 
q^{p_1+\cdots +p_m},
$$
so that by (4.3) we have
$$\displaylines{\quad
{1\over \prod\limits_l{m_l\choose c(B_l)}}
{A'_U(t,q;\bfc)\over (t;q)_{m+1}}=\sum_{s\ge 0} t^s 
\sum_{s\ge p_1\ge \cdots \ge p_m\ge 0} q^{\Sigma\,p_i}
\sum_{\overline w} t^{\des'_U \overline w}
q^{\maj_U'\overline w}\hfill\cr
\kern1.2cm{}
=\sum_{s',{\bf p},\overline w} t^{s'+\des'_U \overline w}
q^{\Sigma\,p_i+\maj_U'\overline w}
=\sum_{(s,(a_{l,i}))} t^s 
q^{\Sigma a_{l,i}}\hfill[{\rm by\ (4.5)}]\cr
\kern1.2cm{}
=\sum_{s\ge0} t^s
\prod_{l;\beta_l=0} 
\sum_{s\ge a_{l,1}\ge\cdots\ge a_{l,m_l}\ge 0}
q^{a_{l,1}+\cdots+a_{l,m_l}} \hfill\cr
\kern4cm{}
\times
\prod_{l;\beta_l=1} 
\sum_{s\ge a_{l,1}>\cdots>a_{l,m_l}\ge 0}
q^{a_{l,1}+\cdots+a_{l,m_l}}\hfill\cr
\kern1.2cm{}
=\sum_{s\ge0} t^s
\prod_{l;\beta_l=0}  {m_l+s\brack m_l}
\prod_{l;\beta_l=1} q^{m_l\choose 2} {s+1\brack m_l}\hfill
[{\rm by\ (2.6)\ and\ (2.7)}.]\cr
}
$$
Hence (4.1) is established.\halmos

\goodbreak

\medskip
As (4.2) implies (2.10) and as the latter identity holds in the
$U$-number-of-inversion version, we have another proof of the
`if' part of Theorem~1.

\bigskip
\centerline{\bf 5. The bijective proof of Theorem 1}

\medskip
Let $U$ be a bipartitional relation.
In this section we construct a bijection $\Phi_U$ of each class
$R(\bfc)$ onto itself satisfying 
$$
\maj'_U w=\inv'_U \Phi_U(w). \leqno(5.1)
$$
One of the main ingredients in the construction of $\Phi_U$ is
the second fundamental transformation~$\Phi$ (see, e.g., [Lo83,
chap.~10]) that satisfies
$$
\maj w=\inv \Phi(w) \leqno(5.2)
$$
on each rearrangement class. The bijection $\Phi_U$ is the
{\it conjugate} of $\Phi$ in the sense that we have
$$
\Phi_U=\delta^{-1}\circ\Phi\circ\delta, \leqno(5.3)
$$
for a certain bijection $\delta$.

Let us first recall the construction of $\Phi$ [Lo83,
chap.~10]: let $w$ be a word in the alphabet~$X$ and $x\in X$.
Two cases are to be considered

\quad(i) the last letter of $w$ is greater than $x$;

\quad(ii) the last letter of $w$ is at most equal to $x$.

\noindent
In case (i) let $(w_1x_1,w_2x_2,\ldots, w_hx_h)$ be the
factorization of~$w$ having the following properties: $x_1$,
$x_2$, \dots~, $x_h$ are letters of~$X$ {\it greater} than~$x$
and $w_1$, $w_2$, \dots~, $w_h$ are words all letters of which
are less than or equal to~$x$.

In case (ii) $x_1$, $x_2$, \dots~, $x_h$ are letters of~$X$ at
most equal to~$x$, while $w_1$, $w_2$, \dots~, $w_h$ are words
all letters of which are greater than~$x$. 

Call $x$-{\it factorization} the above factorization. In both
cases we have
$$
\eqalignno{w&=w_1x_1w_2x_2\ldots w_hx_h\,;\cr
\noalign{\hbox{then define}}
\gamma_x \,w&=x_1w_1x_2w_2\ldots x_hw_h.\cr} 
$$
The construction of $\Phi$ goes as follows. If $w$ is of
length~0 or~1, let $\Phi(w)=w$. For a word $wx$ with $x\in X$
and $w$ of positive length, form $\Phi(w)$ (already defined by
induction), apply $\gamma_x$ to $\Phi(w)$ and add $x$ at the
end of the resulting words, i.e., define
$$
\Phi(wx)=\bigl(\gamma_x\,\Phi(w)\bigr)x. 
$$

Property (5.2) was proved in [Fo68] (see also [Lo83, chap.~10]).
We shall make use of two further properties.

\proclaim Proposition 5.1. 

{\rm (i)} Both $w$ and $\Phi(w)$ end with the
same letter.

{\rm (ii)} Let $y$ and $y'$ be two adjacent letters (with respect
to the usual order) in the alphabet~$X$ and suppose
that both occur exactly once in~$w$. Then, if $y$
occurs to the left of~$y'$ in~$w$, the same holds for $\Phi(w)$.
\endproclaim

Property (i) is true by the very definition of~$\Phi$. Property
(ii) requires a simple verification that will be left out.\halmos

\medskip
Let $((B_1,\ldots,B_k)(\beta_1,\ldots,\beta_k))$ be the ordered
bipartition corresponding to the bipartitional relation~$U$. We
keep the notations given in \S\kern2pt 2. If $w$ is a
word in $R(\bfc)$, let $m_l=\sum B_l$ be the number of letters
in~$w$ belonging to~$B_l$ and let $w(B_l)$ be the {\it subword}
of~$w$ consisting of all the letters belonging to~$B_l$
$(l=1,\ldots, k)$. 

\medskip
The conjugation $\delta$ is defined as follows.

\quad (i) For every $l=1,\ldots, k$ replace each letter
belonging to~$B_l$ by $b_l=\min B_l$ (with respect to the
usual order). Call $\overline w$ the resulting word.

\quad (ii) If $l$ is {\it non-underlined}, read $\overline w$
{\it from left to right} and replace the successive occurrences
of~$b_l$ by $(b_l,1)$, $(b_l,2)$, \dots~, $(b_l,m_l)$; do this for
each non-underlined~$l$.

\quad (iii) Do the operation described in (ii) for each
{\it underlined}~$l$, but this time read $\overline w$ {\it from
right to left}.

The word derived after all those operations will be denoted by
$w_U$. It is actually a rearrangement of the word
$(b_1,1)\ldots(b_1,{m_1})\ldots(b_k,1)\ldots(b_k,m_k)$ (all
letters distinct.) Furthermore, $w_U$ contains the subword 
$(b_l,1)\ldots(b_l,{m_l})$ (resp. $(b_l,m_l)\ldots (b_l,1)$) if $l$
is non-underlined (resp. underlined). To be able to define ``maj"
for $w_U$ we need a {\it linear} order on those ordered pairs. We
shall take:
$$
(b_l,j)>(b_{l'},j')\ \hbox{iff
$l<l'$ or $l=l'$ and $j>j'$.}\leqno(5.3)
$$

The conjugation $\delta$ is then defined by
$$
\delta: w\mapsto (w_U, w(B_1),w(B_2),\ldots,w(B_r)). \leqno(5.4)
$$
The inverse map $\delta^{-1}$ simply consists of replacing
each {\it subword}
$$(b_l,1)(b_l,2)\ldots(b_l,{m_l})\qquad{\rm (resp.}\ 
(b_l,m_l)\ldots  (b_l,2)(b_l,1))
$$
within $w_U$ by the subword $w(B_l)$.

\proclaim Lemma 5.2. With ``$\maj$" defined by
means of the total order $(5.3)$ the following identity holds:
$$
\maj'_U w=\maj w_U.\leqno(5.5)
$$
\endproclaim

\proof
Let $w=x_1x_2\ldots x_m$ and $w_U=z_1z_2\ldots z_m$ (the
letters $z_i$ are ordered pairs $(b_l,j)$). If $(x_i,x_{i+1})\in
U$, then either $x_i\in B_l$, $x_{i+1}\in B_{l'}$ with $l<l'$,
or $x_i$, $x_{i+1}$ are both in the same {\it underlined}
block~$B_l$.

In the first case, $z_i=(b_l,j)$ and $z_{i+1}=(b_{l'},j')$ for
some $j,j'$. But as $l<l'$, we have $z_i>z_{i+1}$ by (5.3).
In the second case, $z_i=(b_l,j)$ and $z_{i+1}=(b_l,j')$; but as
$l$ is underlined we have $j>j'=j-1$ and again
$z_i>z_{i+1}$.

Now if $(x,x_{i+1})\not\in U$, then either 
$x_i\in B_l$, $x_{i+1}\in B_{l'}$ with $l>l'$,
or $x_i$, $x_{i+1}$ are both in the same {\it non-underlined}
block~$B_l$. In the first case the same argument as above
shows that $z_i<z_{i+1}$. In the second case the labelling
from left to right of the non-underlined letters of $\overline w$
yields $z_i=(b_l,j)<(b_l,j+1)=z_{i+1}$.\halmos

\medskip
Next apply the second fundamental transformation to $w_U$. We
obtain a rearrangement $\Phi(w_U)$ that statisfies 
$$
\maj w_U=\inv\Phi(w_U).\leqno(5.6)
$$

\proclaim Lemma 5.3. For each $l=1,2,\ldots, k$ both words
$w_U$ and $\Phi(w_U)$ contain the subword
$$
(b_l,1)(b_l,2)\ldots(b_l,{m_l})\qquad{\rm (resp.}\ 
(b_l,m_l)\ldots  (b_l,2)(b_l,1))
$$
depending on whether $l$ is non-underlined or underlined.
\endproclaim

\proof
This is a consequence of Proposition 5.1 (ii).\halmos

\medskip
Finally, if we apply the conjugation $\delta^{-1}$ to
$\Phi(w_U)$ using the subwords $w(B_1)$, \dots~, $w(B_k)$, we
obtain a rearrangement $\delta^{-1}\,\Phi(w_U)$ 
which is a rearrangement of the original word~$w$ and
satisfies
$$
\inv \Phi(w_U) = \inv_U \delta^{-1}\,\Phi(w_U).\leqno(5.7)
$$
We shall denote it by $\Phi_U(w)$. All the above transformations
are reversible. The product 
$\Phi_U=\delta^{-1}\circ\Phi\circ\delta$ is a
well-defined bijection of $R(\bfc)$ onto itself satisfying (5.1).

\bigskip
\centerline{\bf 6. A proof of the `only if' part of Theorem 1}

\medskip
The proof of that ``only if" part will be the consequence of the
following sequence of lemmas.

\proclaim Lemma 6.1. If there exists an element $x\in X$ such
that  $U\subset (X\setminus \{x\})\times (X\setminus \{x\})$ and
$U\not= \emptyset$, then the equidistribution of $\inv_U$ and
$\maj_U$ does not hold.
\endproclaim

\proof
Let $w$ be a word having no letter equal to~$x$ and let $v$ be a
word in the class $R(x^mw)$ $(m\ge 1)$. Denote by $\overline v$
the word derived from a word~$v$ by deleting all its letters
equal to~$x$. Then $\inv_U v=\inv_U \overline v$. On the other
hand,  $\maj_U x^m\overline v=m\times \des_U \overline
v+\maj_U\overline v$. As $U$ is non-empty, there exists a
rearrangement class $R(w)$ and a word $w'\in R(w)$ such that
$\des_U w'\ge 1$. Thus there is a bound $b$ such that for every
$m\ge 1$ and for every $v\in R(x^mw)$ we have
$\inv_U v\le b$, while
$$
\max_{v\in R(x^mw)} \maj_U v
\ge \maj_U x^m\overline v\ge m.
\qed
$$

\proclaim Lemma 6.2. If $(x,y)\in U$, $(y,x)\in U$ and $x\not=y$,
and if the equidistribution of $\inv_U$ and $\maj_U$ holds, then
$$
(x,x)\in U \quad {\rm and}\quad (y,y)\in U. 
$$
\endproclaim

\proof
Suppose $(x,x)\not\in U$. In the class $R(x^2y)$ we have
$\inv_U w=2$ for all~$w$, while $\maj_U xyx=3$
and the equidistribution does not hold for $R(x^2y)$.\halmos

\medskip
Let $X=\{x,y,z\}$; at this stage it would be useful to have a
thorough table of the relations~$U$ on $X\times X$ for which
the equidistribution of $\inv_U$ and $\maj_U$ holds. As there
are six elements in $X\times X\setminus \diag X\times X$,
there would be only sixty-four cases to consider. As there are
many symmetries, the table could be rapidly set up. A
verification by computer could also be used.
We have preferred to verify the property in each case.

\proclaim Lemma 6.3. If $(x,y)$, $(y,x)$, $(x,z)$ and $(z,x)$ belong
to~$U$, if $x$, $y$ and $z$ are different and if
the equidistribution of $\inv_U$ and $\maj_U$ holds, then $U$
contains the product $\{x,y,z\}\times \{x,y,z\}$.
\endproclaim

\proof
In other words, besides $(x,x)$, $(y,y)$ and $(z,z)$ (as shown in
Lemma~6.2), the relation~$U$ must also contain $(y,z)$ and
$(z,y)$.

If the conclusion does not hold, there are three cases to be
studied. The distributions of $\inv_U$ and $\maj_U$ on the
rearrangement class $R(xyz)$ are shown in the next table and
are never identical.\halmos

\bigskip
\centerline{%
\vbox{\halign{\vrule\strut\ \hfil$#$\hfil\ \vrule
&\kern 13pt \hfil$#$\hfil\kern 13pt \vrule 
&\kern 13pt \hfil$#$\hfil\kern 13pt \vrule 
&\kern 13pt \hfil$#$\hfil\kern 13pt \vrule 
&\kern 13pt \hfil$#$\hfil\kern 13pt \vrule 
&\kern 13pt \hfil$#$\hfil\kern 13pt \vrule 
&\kern 13pt \hfil$#$\hfil\kern 13pt \vrule \cr
\noalign{\hrule}
&\multispan2 \hfil $(y,z)\not\in U, (z,y)\in U$\hfil\ \vrule
&\multispan2 \hfil $(y,z)\in U, (z,y)\not\in U$\hfil\ \vrule
&\multispan2 \hfil $(y,z)\not\in U, (z,y)\not\in U$\hfil\ \vrule\cr
\noalign{\hrule}
w&\inv_U&\maj_U&\inv_U&\maj_U&\inv_U&\maj_U\cr
\noalign{\hrule}
xyz&2&1&3&3&2&1\cr
\noalign{\hrule}
xzy&3&3&2&1&2&1\cr
\noalign{\hrule}
yxz&2&3&3&3&2&3\cr
\noalign{\hrule}
yzx&2&2&3&3&2&3\cr
\noalign{\hrule}
zxy&3&3&2&3&2&3\cr
\noalign{\hrule}
zyx&3&3&2&2&2&2\cr
\noalign{\hrule}
}}}

\bigskip
If $U$ is a relation on $X\times X$, its {\it symmetric} part,
i.e., the set of all ordered pairs $(x,y)$ such that both $(x,y)$
and $(y,x)$ belong to~$U$, is denoted by $S(U)$. Also let
$A(U)=U\setminus S(U)$ be its {\it asymmetric part}. Finally, let
$X_U$ be the subset of~$X$ of all the $x$'s such that $(x,y)\in
S(U)$ (and so $(y,x)\in S(U)$) for some $y\in X$.

\proclaim Lemma 6.4. If the equidistribution holds for $U$, then
$S(U)$ is an equivalence relation on $X_U\times X_U$.
\endproclaim

\proof
Let $x\in X_U$ and let $y\in X$ such that $(x,y)\in S(U)$. If
$y=x$, then $(x,x)\in U$. If $y\not=x$, Lemma~6.2 also implies
that $(x,x)\in U$. Thus $S(U)$ is reflexive. By definition, $S(U)$
is symmetric. Now let $x,y,z\in X_U$ and suppose $(x,y)\in S(U)$
and $(y,z)\in S(U)$. Then Lemma~6.3 implies that $(x,y)\in S(U)$.
The relation is then transitive.\halmos

\medskip
Thus, if the equidistribution holds for $U$, there is a partition 
$\{B_1,\ldots,B_l\}$ of $X_U$ such that
$S(U)=B_1\times B_1 \cup \cdots \cup B_l\times B_l$.
The subsets $B_1$, \dots~, $B_l$ will be called the {\it blocks
of}~$X_U$.

\proclaim Lemma 6.5. Suppose that the equidistribution holds for
$U$ and let $x,y$ be two distinct elements belonging to the same
block, say, $B_i$ of $X_U$ and let $z$ be an element of~$X$. Then
$$
\displaylines{
(z,x)\in U\Leftrightarrow (z,y)\in U\,;\cr
(x,z)\in U\Leftrightarrow (y,z)\in U.\cr} 
$$
\endproclaim

\proof
If $z\in B_i$, then $(x,z)\in S(U)$ and $(y,z)\in S(U)$ and there
is nothing to prove. If $z$ belongs to another block $B_j$ of
$X_U$ and  if $(z,x)\in U$, then $(x,z)\not\in U$. Otherwise, we
would have $B_i=B_j$. If $z\not\in X_u$ and $(z,x)\in U$, again
$(x,z)\not\in U$. Otherwise, $(x,z)\in S(U)$ and this would
contradict $z\not\in X_U$.

Suppose that conditions $(z,x)\in U$ and $(z,y)\not\in U$
hold. Two cases are to be considered : 
(a) $(y,z)\not\in U$ ;
(b) $(y,z)\in U$.

In case (a) we have $(x,y), (y,x),(z,x)\in U$,
$(x,z),(z,y),(y,z)\not\in U$. But for each word~$w$ in the class
$R(xyz)$ we have $\inv_U w\le 2$, while $\maj_U zxy=3$, so that
the equidistribution would not hold.

In case (b) we have $(x,y), (y,x),(z,x),(y,z)\in U$,
$(x,z),(z,y)\not\in U$. Let $V=X\times X\setminus U$. Then
$\inv_V w\le 2$ for all~$w$, while $\maj_V xzy=3$. Thus, the
equidistribution would not hold for~$V$ and also for~$U$.

Thus cases (a) and (b) cannot occur and consequently if
$(z,y)\in U$ holds, we must have $(z,y)\in U$. The elements $x$
and~$y$ play a symmetric role, so that the first equivalence is
proved.

The proof of the second equivalence is quite analogous. If
$(x,z)\in U$, we have seen that $(z,x)\not\in U$.
Suppose $(y,z)\not\in U$ and consider the two cases :
(a) $(z,y)\not\in U$ ; (b) $(z,y)\in U$.

In case (a) we have $(x,y),(y,x),(x,z)\in U$,
$(z,x),(y,z),(z,y)\not\in U$. Again for each $w\in R(xyz)$ we have
$\inv_U w\le 2$, while $\maj_U yxz=3$, so that the
equidistribution does not hold.

In case (b) we have $(x,y),(y,x),(x,z),(z,y)\in U$,
$(z,x),(y,z)\not\in U$. Let $V=X\times X\setminus U$. 
Then $\inv_V w\le 2$, while $\maj_V yzx=3$, so that the
equidistribution does not hold for~$V$, and then for~$U$. 

As $x$ and $y$ play a symmetric role, the second equivalence is
also established.\halmos

\medskip
{\bf Notation}. It will be convenient to write $x\rightarrow y$
for $(x,y)\in U$, $x\not\rightarrow y$ or
$y\nonfleche x$ for $(x,y)\not\in U$,
$x\rightleftharpoons y$ for $(x,y)\in S(U)$.

\proclaim Lemma 6.6. Suppose that the equidistribution holds
for~$U$. Then, either there is a block $B_1$ of $X_U$ with the
property

$\forall\,x\in B_1$, $\forall\,y\in X\setminus B_1$, then
$(x,y)\not\in U$ ;

\noindent
or there exists an $x\in X\setminus X_U$ such that

$\forall\,y\in X$, then $(x,y)\not\in U$.
\endproclaim

The foregoing property means that by rearranging the elements
of~$X$, either the top left corner $B_1\times (X\setminus B_1)$
of $X\times X$, or a left block $C\times X$ has no intersection
with~$U$.

\medskip
\proof
Suppose that the conclusion is false. This means that for every
block $B_i$ of $X_U$ there is $x\in B_i$, $y\in X\setminus
B_i$ such that $(x,y)\in U$ and also that for all $x\in
X\setminus X_U$ there is $y\in X\setminus\{x\}$ such that
$(x,y)\in U$.

Let $x_0\in X$. If $x_0$ belongs to a block $B_{i_0}$, there is
$x_1\in B_{i_0}$ and $x_2\not\in B_{i_0}$ such that
$x_0\rightleftharpoons x_1\rightarrow x_2$. But the previous
lemma says that :

\bigskip
\newbox\partitionalC

\setbox\partitionalC=\vbox{\offinterlineskip 
\segment(0,0)\dir(0,-1)\long{30}
\segment(0,0)\dir(1,0)\long{30}
\segment(30,0)\dir(0,-1)\long{30}
\segment(0,-30)\dir(1,0)\long{30}
\segment(10,0)\dir(0,-1)\long{20}
\segment(10,-20)\dir(1,0)\long{20}
\segment(10,-30)\dir(0,1)\long{10}
\centerput(20,-25){$(X\setminus C)\times C$}
\centerput(5,-34.5){$C$}
\centerput(20,-34.5){$X\setminus C$}
\centerput(5,-15){$C\!\!\times \!\!X$}
}

\newbox\partitionalB

\setbox\partitionalB=\vbox{\offinterlineskip 
\segment(0,0)\dir(0,-1)\long{30}
\segment(0,0)\dir(1,0)\long{30}
\segment(30,0)\dir(0,-1)\long{30}
\segment(0,-30)\dir(1,0)\long{30}
\segment(10,0)\dir(0,-1)\long{20}
\segment(0,-20)\dir(1,0)\long{30}
\centerput(15,-25){$X\times B_1$}
\centerput(5,-34.5){$B_1$}
\centerput(20,-34.5){$X\setminus B_1$}
\centerput(5,2){$B_1\!\!\times\!\! (X\!\setminus\! B_1)$}
}

\centerline{\box\partitionalB\hskip40mm\box\partitionalC\hskip
30mm} 
\kern 105pt
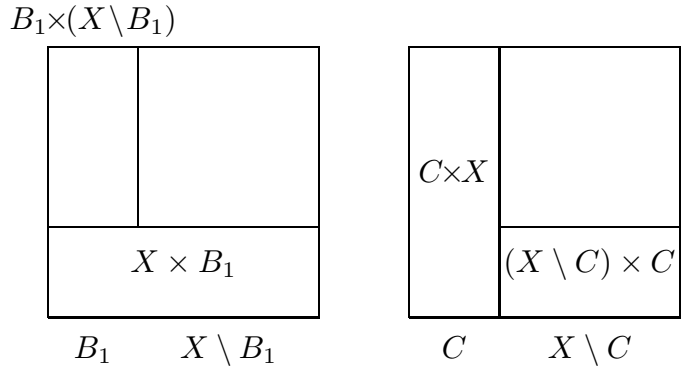
\centerline{Fig. 2}

\bigskip
if $x_0\rightleftharpoons x_1$ and $x_1\rightarrow x_2$,
then $x_0\rightarrow x_2$. Also $x_2\not\rightarrow x_0$.

\noindent
Now, either $x_2\in B_{i_2}$ with $i_2\not=i_1$ or $x_2\not\in
X_U$. In the first case there is $x_3\in B_{i_2}$ and
$x_4\not\in B_{i_2}$ such that
$x_2\rightleftharpoons x_3\rightarrow x_4$.
Using the same lemma we also have
$x_2\rightarrow x_4$ and $x_4\not\rightarrow x_2$.
If $x_2\not\in X_U$, there is $x_4\not=x_2$ such that
$x_2\rightarrow x_4$ and also $x_4\not\rightarrow x_2$, because
$x_2\not\in X_U$.

We can then build a sequence $(x_0,x_2,x_4,\ldots\,)$ with the
property
$$
\displaylines{
x_0\rightarrow x_2\rightarrow x_4\rightarrow
x_6\rightarrow\cdots\cr
x_0\nonfleche  x_2\nonfleche 
x_4\nonfleche  x_6\nonfleche\cdots\cr}
$$
and such that $x_{2i}\not=x_{2i+2}$ at each step. If we had
started with an element $x_0\notin X_U$, the conclusion would
have been the same.

The above sequence cannot be infinite and have all its elements
distinct, so that, after relabelling, there is a finite sequence
$(y_1,y_2,\ldots,y_{n+1})$ of elements of~$X$ with the following
properties :

(a) $n\ge 2$ ;

(b) all terms $y_1$, $y_2$, \dots~, $y_n$ are different ;

(c) $y_1\rightarrow y_2\rightarrow y_3\rightarrow
\cdots \rightarrow y_n\rightarrow y_{n+1}=y_1$ ;

(d) $y_1\nonfleche y_2\nonfleche y_3\nonfleche
\cdots \nonfleche y_n\nonfleche y_{n+1}=y_1$ ;

\noindent
If $n=2$, we have $y_1\rightarrow y_2\rightarrow y_1$
and $y_1\nonfleche y_2\nonfleche y_1$, a contradiction, so that
$n\ge 3$.

Consider the class $R(y_1y_2\ldots y_n)$. Then
$\maj_U y_1y_2\ldots y_n=1+2+\cdots+(n-1)=n(n-1)/2$. 
If there is a word $v=z_1z_2\ldots z_n\in R(y_1y_2\ldots y_n)$
such that $\inv_U v=n(n-1)/2$, this means that
$$\eqalign{
z_1\rightarrow z_2,\ z_1\rightarrow z_3,\ \ldots ,
z_1&\rightarrow z_n\cr
z_2\rightarrow z_3,\ \ldots ,
z_2&\rightarrow z_n\cr
\ldots&\cr
z_{n-1}&\rightarrow z_n\cr
}\eqno{(\star)}
$$
If $z_1=y_i$ with $2\le i\le n+1$, let $j$ be
the unique integer such that $z_j=y_{i-1}$.
Relation~(e) above says that
$z_j=y_{i-1}\nonfleche y_i=z_1$. This contradicts
$(\star)$.\halmos

Let $C=\{x\in X\setminus X_U : \forall\,y\in X,\
x\not\rightarrow y\,\}$. From Lemma~6.6 it follows that, if the
equidistribution hods and $C$ is empty, there is a unique block
$B_1$ of $X_U$ such that
$$
\forall\, x\in B_1,\ \forall\, y\in X\setminus B_1,\quad
{\rm then}\ x\not\rightarrow y.\leqno(\star\,\star) 
$$

\proclaim Lemma 6.7. Suppose that the equidistribution holds
for~$U$. If $C$ is non-empty, then
$$
(X\setminus C)\times C\subset U. \leqno{(\star\,\star\,\star)}
$$
In other words, $\forall\, y\in X\setminus C$, $\forall\, x\in C$,
then $y\rightarrow x$.

If $C$ is empty and if $B_1$ is the block defined in
$(\star\,\star)$, then
$$
(X\setminus B)\times B \subset U. \leqno{(iv)}
$$
In other words, $\forall\, y\in X\setminus B$, $\forall\, x\in B$,
then $y\rightarrow x$.
\endproclaim

See Fig. 1 : in each case the bottom rectangle to the right is
entirely contained in~$U$.

\medskip
\proof
Assume that $C$ is non-empty and suppose that
$(\star\,\star\,\star)$ does not hold. Then there is $y\in
X\setminus C$ and also $x\in C$ such that $y\not\rightarrow x$
and $y\not= x$. As $y\not\in C$, there is $z\in X$ such that
$y\rightarrow z$. Notice that $z$ may be equal to $y$, if $y\in
X_U$, but $z\not=x$, as we have assumed $y\not\rightarrow x$.

\newbox\partitionalA

\setbox\partitionalA=\vbox{\offinterlineskip 
\segment(0,0)\dir(0,-1)\long{30}
\segment(0,0)\dir(1,0)\long{30}
\segment(30,0)\dir(0,-1)\long{30}
\segment(0,-30)\dir(1,0)\long{30}
\segment(0,-20)\dir(1,0)\long{30}
\segment(0,-25)\dir(1,0)\long{30}
\segment(5,-30)\dir(0,1)\long{5}
\segment(10,-30)\dir(0,1)\long{10}
\segment(20,-30)\dir(0,1)\long{5}
\segment(25,-30)\dir(0,1)\long{5}
\centerput(5,-32.5){$x$}
\centerput(20,-32.5){$y$}
\centerput(25,-32.5){$z$}
\centerput(-3,-25){$x$}
\centerput(5,-35.5){$C$ or $B_1$}
}

\centerline{\box\partitionalA\hskip30mm}

\kern 105pt 

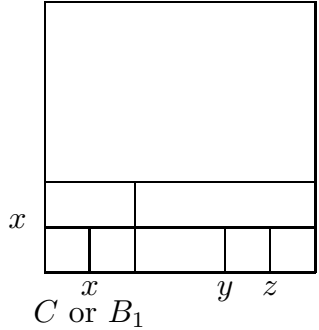
\centerline{Fig. 3}

\smallskip
Consider the class $R(xyz)$. By assumption, $y\not\rightarrow x$,
$y\rightarrow z$ and also $x\not\rightarrow y$,
$x\not\rightarrow z$, since $x\in C$. Four cases are to be
considered :

(a) $z\rightarrow x$, $z\rightarrow y$ ;
(b) $z\rightarrow x$, $z\not\rightarrow y$ ;
(c) $z\not\rightarrow x$, $z\rightarrow y$ ;
(d) $z\not\rightarrow x$, $z\not\rightarrow y$.

In both cases (a) and (b) $\inv_U w\le 2$ for all~$w$, while
$\maj_U yzx=3$. In both cases (c) and (d) $\inv_U w\le 1$ for
all~$w$, while $\maj_U xyz=2$. Thus there is never
equidistribution on $R(xyz)$.

Suppose that $C$ is empty. Let $B_1$ be the block defined in
$(\star\,\star)$. If $(iv)$ does not hold, there is $y\in
X\setminus B_1$ and also $x\in B_1$ such that
$y\not\rightarrow z$. As $y\notin B_1$ and since $C$ is
supposed to be empty, there exists $z$ such that 
$y\rightarrow z$. Again we have
$y\not\rightarrow x$, $y\rightarrow z$, $x\not\rightarrow y$,
$x\not\rightarrow z$. The same analysis as above shows that
there is no equidistribution on $R(xyz)$.\halmos

It follows from Lemma 6.6 and Lemma 6.7 that, if the
equidistribution holds for~$U$, then, either
$C$ is non-empty and then

$C\times X$ is empty and $(X\setminus C)\times C\subset U$,

\noindent
or $C$ is empty and then there is a unique block $B_1$
of $X_U$ such that

$B_1\times (X\setminus B_1)$ is empty and
$(X\setminus B_1)\times \subset U$.

The theorem is now easily proved by induction on
$\card X$. If the equidistribution holds for~$U$ defined on
$X\times X$ and if $C$ is non-empty, then the equidistribution
also holds for the relation $V=U\cap(X\setminus C)\times
(X\setminus C)$ defined on $(X\setminus C)\times
(X\setminus C)$. By induction $V$ is bipartitional. Hence, $U$ is
also bipartitional (see Fig.~1).

In the same manner, if $C$ is empty, then the equidistribution
holds for the relation $V=U\cap(X\setminus B_1)\times
(X\setminus B_1)$. By induction $V$ is bipartitional. Hence, $U$ is
also bipartitional.

\goodbreak
\bigskip
\centerline{\bf 7. Compatible bipartitional relations}

\medskip
The statistics ``$\maj_U$" and ``$\inv_U$" have been defined in
(0.2); also remember that ``$\des_U$" counts the $U$-descents
of the second kind, as defined at the end of the introduction. The
calculation of the generating function for $(\des_U,\maj_U)$ and
the construction of the bijection that carries ``$\maj_U$" onto
``$\inv_U$" will be very similar to their equivalent derivations
for $(\des'_U,\maj'_U)$, ``$\maj'_U$" and ``$\inv_U'$." Let
$$
\leqalignno{A_U(q;\bfc)&=\sum_w q^{\inv_U w}
\qquad (w\in R(\bfc))\,;&(7.1)\cr
A_U(t,q;\bfc)&=\sum_w t^{\des_U w} q^{\maj_U w}
\qquad (w\in R(\bfc)).&(7.2)\cr}
$$
The identity
$$
\displaylines{\rlap{\rm (7.3)}\hfill
 A_U(q;{\bf c})={|\bfc| \brack m_1,\ldots, m_k}
\prod_{l=1}^k {m_l \choose c(B_l)} 
q^{\beta_l{m_l+1\choose 2}};\hfill\cr
}
$$
follows from (2.10), as we have to add the
total number of underlined letters, i.e.,
$\sum\limits_l \beta_l m_l$ to the power of~$q$.

The proof of the formula
$$\displaylines{
 \rlap{\rm (7.4)}\hfill
\eqalign{
\sum_{\bf c} {A_U(q;{\bf c})\over
(q)_{\left |{\bf c}\right|}} {\bf u}^{\bf c}
&=\prod_{l;\beta_l=0} e_q(\scriptstyle\sum \textstyle u(B_l))
\times \displaystyle \prod_{l;\beta_l=1}
E_q(q\scriptstyle\sum \textstyle u(B_l))\hfill\cr 
&={\prod\limits_{l;\beta_l=1} (-q\sum u(B_l);q)_\infty
\over \prod\limits_{l;\beta_l=0} (\sum u(B_l);q)_\infty}\,;\cr
}\hfill\cr
} 
$$
follows the same pattern as the proof of (2.11).

Let $A^{\maj}_U(q;\bfc)=\sum\limits_w q^{\maj_U w}$ $(w\in
R(\bfc))$. Again, we don't prove that $A^{\maj}_U(q;\bfc)$ is
equal to the right-hand side of (7.3). We'd rather derive the
formulas for $A_U(t,q;\bfc)$, defined in (7.2), in the spirit of 
section~4.

\proclaim Proposition 7.1. Let $U$ be a compatible
bipartitional relation. Then
$$\displaylines{(7.5)\quad
{A_U(t,q;\bfc)\over (t;q)_{|\bfc|+1}}\hfill\cr
\hfill{}
=\prod_{l=1}^k 
{ m_l\choose c(B_l)}
\times \sum_{s\ge 0} t^s \prod_{l;\beta_l=0}
{m_l+s\brack m_l}\times
\prod_{l;\beta_l=1} q^{m_l+1\choose 2} {s\brack m_l};\quad
\cr
\rlap{\rm (7.6)}\hfill
\sum_\bfc {A_U(t,q;\bfc)\over (t;q)_{|\bfc|+1}} {\bf u}^\bfc
=\sum_{s\ge 0} t^s\, {\prod\limits_{l;\beta_l=1}
(-q\sum u(B_l);q)_{s}\over
\prod\limits_{l;\beta_l=0}
(\sum u(B_l);q)_{s+1}}.\hfill\cr}
$$
\endproclaim

\proof
Let $w\mapsto \bigl(\overline w, w(B_1),\ldots,w(B_l)\bigr)
\quad{\rm and}\quad 
(s',{\bf p},\overline w)\mapsto (s,(a_{l,i}))$
be the two bijections defined in section~4. We keep the same
notations as in that section. In particular, let
$\overline w= \overline x_1\overline x_2 \ldots \overline x_m$.
The only difference to be brought to the constructions of those
bijections is to notice that $z_m=1$ iff $\overline x_m$ belongs
to an underlined block. Consequently, the sequences
$a_{l,1}\ldots a_{l,m_l}$ associated with the {\it underlined}
blocks are still {\it strictly} decreasing, but also $a_{l,m_l}\ge
1$.

The reason is the following: let $l$ be underlined and
let $\overline x_i$ be the {\it rightmost} letter of~$\overline
w$ that belongs to the block~$B_l$. If $i=m$, then
$a_{l,m_l}=p_m+z_m\ge 1$; if $i<m$, then, either there is one
{\it non-underlined} letter in the factor $\overline
x_{i+1}\ldots\overline x_m$ and necessarily one
$U$-descent because $U$ is supposed to be {\it compatible}, or
all the letters in that factors are underlined and in particular
$z_m=1$. In both cases, $a_{l,m_l}\ge 1$.

Accordingly, the mapping $(s',{\bf p},\overline w)\mapsto
(s,(a_{l,i}))$ is a bijection satisfying
$$\displaylines{
\eqalign{s&\ge a_{l,1}\ge \cdots \ge a_{l,m_l}\ge 0,
\hbox{ if $l$ is {\it not} underlined;}\cr
s&\ge a_{l,1}>\cdots >a_{l,m_l}\ge 1,
\hbox{ if $l$ is underlined.;}\cr
}\cr
s=s'+\des_U \overline w\,;\cr
\sum_{l,i}a_{l,i}=
\sum_i p_i +\maj_U \overline w.\cr}
$$
In the same manner as in section 4 we have
$$\displaylines{\quad
{1\over \prod\limits_l{m_l\choose c(B_l)}}
{A_U(t,q;\bfc)\over (t;q)_{m+1}}=\sum_{s\ge 0} t^s 
\sum_{s\ge p_1\ge \cdots \ge p_m\ge 0} q^{\Sigma\,p_i}
\sum_{\overline w} t^{\des_U \overline w}
q^{\maj_U\overline w}\hfill\cr
\kern1.5cm{}
=\sum_{s\ge0} t^s
\prod_{l;\beta_l=0} 
\sum_{s\ge a_{l,1}\ge\cdots\ge a_{l,m_l}\ge 0}
q^{a_{l,1}+\cdots+a_{l,m_l}} \hfill\cr
\kern4cm{}
\times
\prod_{l;\beta_l=1} 
\sum_{s\ge a_{l,1}>\cdots>a_{l,m_l}\ge 1}
q^{a_{l,1}+\cdots+a_{l,m_l}}\hfill\cr
\kern1.5cm{}
=\sum_{s\ge0} t^s
\prod_{l;\beta_l=0}  {m_l+s\brack m_l}
\prod_{l;\beta_l=1} q^{m_l+1\choose 2} {s\brack m_l},\hfill
\cr
}
$$
by (2.6) and (2.7).\halmos

\medskip
As (7.4) holds (in the $U$-number-of-inversion version) and
since (7.6) implies (7.4), we then have a proof of the `if' part of
Theorem~2.

The proof of the `only if' part is straightforward. Suppose that
$U$ is non-compatible, so that there is an underlined block
$B_l$ to the left of a non-underlined one $B_{l'}$, i.e., $l<l'$.
Take two integers $x\in B_l$ and $x_l'\in B_{l'}$ and consider
the class $R(xx')$ of the two words $xx'$ and $x'x$. 
Then $\inv_U xx'=2$, $\inv_U x'x=1$, while $\maj_U xx'=3$,
$\maj_U x'x=0$.

\bigskip
\centerline{\bf 8. A bijective Proof of Theorem 2}
\medskip
Let $\pi=(\underline {B_1},\ldots,\underline{B_n},
B_{n+1},\ldots,B_r)$ be a compatible ordered bipartition having
exactly $n$ underlined blocks lying in the beginning and let
$U$ be the corresponding compatible bipartitional relation.
As done in the
papers by Steingr\`\i msson [St93] and Clarke and Foata ({\it op.
cit.}), let us introduce an extra letter~$\star$ and form the new
compatible ordered bipartition
$$
\pi^\star=(\underline{B_1},\ldots,
\underline{B_n},\{\star\},\{B_{n+1}\},\ldots,\{B_r\}).\leqno(8.1)
$$
Denote by  $U^\star$ the 
bipartitional relation associated with~$\pi^\star$. Notice
that~$U^\star$ is a relation on $(X\cup\{\star\})\times
(X\cup\{\star\})$. We now make use of the transformation
$\Phi_{U^\star}$ (constructed in section~5) on the words in the
alphabet $X\cup\{\star\}$.

If the word $w=x_1x_2\ldots x_m$ belongs to the class
$R(\bfc)$, form the word $w\star$. Its image under
$\Phi_{U^\star}$ will yield a word of the form $w'\star$, by
Proposition~5.1. There is a
$U^\star$-descent at position~$m$, if and only if $x_m$ is
underlined. Hence 
$$
\maj'_{U^\star } w\star=\maj'_{U} w+m\,\chi(x_m\ {\rm is\
underlined}).\leqno(8.2)
$$
Also adding $\star$ at the end of $w'$ will increase the number
of $U^\star$-inversions by exactly the number of underlined
letters in~$w'$, i.e., $|w'|_-$. Hence
$$
\inv'_{U^\star } w'\star
=\inv'_{U} w'+|w'|_{-}.\leqno(8.3)
$$
Hence
$$\eqalignno{\noalign{\vskip -24pt}
\maj_U w&=\maj'_{U} w+m\,\chi(x_m\ {\rm is\
underlined})&[{\rm by\ (0.2)}]\cr 
&=\maj'_{U^\star } w\star&[{\rm by\ (8.2)}]\cr
&=\inv'_{U^\star } w'\star&[{\rm by}\ (5.1)]\cr
&=\inv'_{U} w'+|w'|_{-}&[{\rm by\ (0.2)}]\cr
&=\inv_{U} w'.&[{\rm by\ (8.3)}]\cr }
$$
As $\Phi_{U^\star}$ maps the set of all words in each
rearrangement class ending with~$\star$ onto the same set, the
mapping $w\mapsto w'$ is a bijection of $R(\bfc)$ onto itself.
Moreover, it satisfies
$$
\maj_U w=\inv_Uw'.\leqno(8.4) 
$$

{\bf Remark:} Formulas (4.2) and (7.6) are the factorial
generating functions for the pairs $(\des'_U,\maj'_U)$ and
$(\des_U,\maj_U)$, respectively. On the other hand, the bijection
$\Phi_U$ (defined in \S\kern2pt 5) and the bijection $w\mapsto
w'$ just defined (that we shall denote by $\Psi_U$)  satisfy
(4.1) and (8.4). Let $\interr'_U=\des'_U\circ \Phi_U^{-1} $ and
$\interr_U=\des_U\circ\Psi_U^{-1} $, so that
$$\eqalign{
(\des'_U,\maj'_U)\,w&=(\interr'_U,\inv'_U)\,\Phi_U(w)\,;\cr
(\des_U,\maj_U)\,w&=(\interr_U,\inv_U)\,\Psi_U(w).\cr
}
$$
Th natural question arises: can we find suitable predicates to
define ``$\interr'_U$" and ``$\interr_U$" without any references
to the bijections $\Phi_U$ and $\Psi_U$?

\def\article#1|#2|#3|#4|#5|#6|#7|
    {{\leftskip=7mm\noindent
     \hangindent=2mm\hangafter=1
     \llap{[#1]\hskip.35em}{#2},\enspace
     ``#3," {\sl #4}, {\bf #5}, {\oldstyle #6},
     p.\nobreak\ #7.\par}}
\def\livre#1|#2|#3|#4|
    {{\leftskip=7mm\noindent
    \hangindent=2mm\hangafter=1
    \llap{[#1]\hskip.35em}{#2},\enspace 
    ``{\sl #3}."\enspace  #4.\par}}
\def\divers#1|#2|#3|
    {{\leftskip=7mm\noindent
    \hangindent=2mm\hangafter=1
     \llap{[#1]\hskip.35em}{#2},\enspace 
     #3.\par}}

\vskip 54pt plus 5pt minus 2pt

\centerline{\bf References}
\vskip 10pt

\livre An76|George E. Andrews|The Theory of
Partitions|London, Addison-Wesley, {\oldstyle 1976}
({\sl Encyclopedia of Math. and Its Appl.}, {\bf 2})|

\livre Be94|Fran\c cois Bergeron, Gilbert Labelle,
Pierre Leroux|Th\'eorie des es\-p\`eces et combinatoire
des structures arborescentes|Montr\'eal, Universit\'e du Qu\'ebec
\`a Montr\'eal, Publ. du LACIM, no.~19, {\oldstyle 1994}|

\article ClFo94|Robert J. Clarke and
Dominique Foata|Eulerian Calculus, I: univariable
statistics|Europ. J. Combinatorics|15|1994|345--362|

\divers ClFo95a|Robert J. Clarke and Dominique
Foata|Eulerian Calculus, II: an extension of Han's fundamental
transformation, to appear in {\sl Europ. J. Combinatorics},
{\oldstyle 1995}|

\divers ClFo95b|Robert J. Clarke and Dominique
Foata|Eulerian Calculus, III: the ubiquitous Cauchy formula,
to appear in {\sl Europ. J. Combinatorics}, {\oldstyle 1995}|

\livre Co70|Louis Comtet|Analyse Combinatoire, {\rm vol. 1
and 2}|Paris, Presses Universitaires de Frances, 
{\oldstyle 1970}. (English edition {\sl Advanced
Combinatorics}, D. Reidel, Dordrecht, {\oldstyle 1974})|

\article Fo68|Dominique Foata|On the Netto inversion 
number of a sequence|Proc. Amer. Math. Soc.|19|1968|236--240|

\livre Fo74|Dominique Foata|La s\'erie g\'en\'eratrice
exponentielle dans les pro\-bl\`emes
d'\'enum\'eration|Montr\'eal, Presses Universitaires de
Mont\-r\'eal, $\oldstyle1974$|

\divers Fo95|Dominique Foata|Les distributions
Euler-Mahoniennes sur les mots, to appear in {\sl Discrete
Math.}, $\oldstyle 1995$|

\livre GaRa90|George Gasper and
Mizan Rahman|Basic Hypergeometric Series|\hfil\break London,
Cambridge Univ. Press, {\oldstyle 1990}  ({\sl Encyclopedia of
Math. and Its Appl.}, {\bf 35})|

\divers Han95|Guo-Niu Han|Ordres bipartitionnaires et
statistiques sur les mots, to appear in {\sl Electronic J.
Combinatorics}, {\oldstyle 1995}|

\livre Kn73|Donald E. Knuth|The Art of Computer
Programming,  {\rm vol.~3}, Sorting and
Searching|Addison-Wesley, Reading,  {\oldstyle1973}|

\article Kn92|Donald E. Knuth|Convolutional
Polynomials|Mathematica J.|2|1992|67--78|

\livre Lo83|M. Lothaire|Combinatorics on words|Reading,
Addison-Wesley, $\oldstyle 1983$  ({\sl Encyclopedia of Math.
and its Appl.}, {\bf 17})|

\article Mac13|P.A. MacMahon|The indices of
permutations and the derivation therefrom of functions of
a single variable associated with the permutations of any
assemblage of objects|Amer. J. Math.|35| 1913|314--321|

\livre Mac15|P.A. MacMahon|Combinatory
Analysis, {\rm vol.~1}|Cam\-bridge, Cambridge Univ. Press,
$\oldstyle 1915$ (Reprinted by Chelsea, New York,
$\oldstyle 1955$)|

\divers Mac78|P.A. MacMahon|{\sl Collected
Papers}, vol.~1 [G.E. Andrews, ed.].\enspace 
Cam\-bridge, Mass., The M.I.T. Press, $\oldstyle 1978$|

\livre Sl94|N.J.A. Sloane|A Handbook of Integer Sequences|New
York, Academic Press, {\oldstyle 1973}|

\article St93|Einar
Steingr\`\i msson|Permutation Statistics of Indexed 
Per\-mu\-ta\-tions|Europ. J. 
Combinatorics|15|1994|xxx--xxx|

\bye